\documentclass[11pt,reqno]{amsart}
\usepackage{amssymb}
\usepackage[all]{xy}
\newtheorem{thm}{Theorem}[section]
\newtheorem{lem}[thm]{Lemma}
\newtheorem{prop}[thm]{Proposition}
\newtheorem{cor}[thm]{Corollary}

\theoremstyle{definition}
\newtheorem{dfn}[thm]{Definition}
\newtheorem{ex}[thm]{Example}

\theoremstyle{remark}
\newtheorem{remark}[thm]{Remark}
\newtheorem{remarks}[thm]{Remarks}

\newcommand{\CA}{{\mathcal{A}}}

\newcommand{\CE}{{\mathcal{E}}}

\newcommand{\CF}{{\mathcal{F}}}

\newcommand{\CK}{{\mathcal{K}}}
\newcommand{\CL}{{\mathcal{L}}}
\newcommand{\CB}{{\mathcal{B}}}
\newcommand{\CG}{{\mathcal{G}}}
\newcommand{\CO}{{\mathcal{O}}}

\newcommand{\CS}{{\mathcal{S}}}

\newcommand{\af}{\alpha}
\newcommand{\bt}{\beta}
\newcommand{\gm}{\gamma}
\newcommand{\dt}{\delta}
\newcommand{\ep}{\varepsilon}

\newcommand{\ld}{\lambda}

\newcommand{\sm}{\sigma}

\newcommand{\om}{\omega}

\newcommand{\Ld}{\Lambda}

\newcommand{\Om}{\Omega}

\begin{document}


\title[Purely infinite labeled graph $C^*$-algebras]
{Purely infinite labeled graph $C^*$-algebras}

\author[J. A. Jeong]{Ja A Jeong$^{\dagger}$}
\thanks{Research partially supported by NRF-2015R1C1A2A01052516$^{\dagger}$}
\thanks{Research partially supported by Hanshin University$^{\ddagger}$}
\address{
Department of Mathematical Sciences and Research Institute of Mathematics\\
Seoul National University\\
Seoul, 08826\\
Korea} \email{jajeong\-@\-snu.\-ac.\-kr }

\author[E. J. Kang]{Eun Ji Kang$^{\dagger}$}
\address{
BK21 Plus Mathematical Sciences Division\\
Seoul National University\\
Seoul, 08826\\
Korea} \email{kkang33\-@\-snu.\-ac.\-kr }

\author[G. H. Park]{Gi Hyun Park$^{\ddagger}$}
\address{
Department of Financial Mathematics\\
Hanshin University\\
Osan, 18101\\
Korea} \email{ghpark\-@\-hs.\-ac.\-kr }




\begin{abstract} 
In this paper, we consider pure infiniteness of 
generalized Cuntz-Krieger algebras 
associated to labeled spaces $(E,\CL,\CE)$. 
It is shown  that 
a $C^*$-algebra $C^*(E,\CL,\CE)$ is purely infinite  in the sense that 
every nonzero hereditary subalgebra contains an infinite projection 
(we call this property (IH)) 
if $(E, \CL,\CE)$ is disagreeable and every vertex connects to a loop.   
We also prove that under the condition  
analogous to (K) for  usual graphs, 
$C^*(E,\CL,\CE)=C^*(p_A, s_a)$ is purely infinite 
in the sense of Kirchberg and R\o rdam 
if and only if every generating projection $p_A$, $A\in \CE$,  is 
properly infinite, and also if and only if  
every quotient of $C^*(E,\CL,\CE)$ 
has the property (IH).
\end{abstract}

\maketitle

\setcounter{equation}{0}

\section{Introduction}

\noindent
In the literature, one can find at least two definitions
for pure infiniteness of $C^*$-algebras. 
The property of a $C^*$-algebra having infinite projections in 
each of its nonzero hereditary subalgebras was called 
`purely infinite' in some articles 
(for example, \cite{A, BPRS, DT, KPR, LS}), 
then the same terminology {\it `purely infinite'} began to be used 
by Kirchberg and R\o rdam \cite{KR} 
to mean another infiniteness of $C^*$-algebras
 wth the motivation  to extend 
to non-simple $C^*$-algebras the fact that 
a simple separable nuclear unital $C^*$-algebra $A$ is 
purely infinite if and only if it is isomorphic to 
$A\otimes \CO_\infty$
  by finding 
the right definition of being purely infinite 
for non-simple $C^*$-algebras.   
While these two definitions are known to be 
equivalent for simple $C^*$-algebras,  
one of the two  is neither 
stronger nor weaker than the other in general, and thus 
we will say to avoid any confusion 
that a $C^*$-algebra has the property (IH)  
if it has only infinite hereditary subalgebras, that is,  
every nonzero hereditary subalgebra contains an infinite projection. 

Pure infiniteness and the property (IH)  have been studied 
for various classes of $C^*$-algebras   
(a few examples are \cite{A, BCS, BOSS, BP2,  EP, H, HZ,  LS, M, RS}
 among many others), 
and especially interesting results are known for graph $C^*$-algebras  
(see \cite{BPRS, DT, KPR}). 
The purpose of the present paper is to do this sort of study 
for labeled graph $C^*$-algebras $C^*(E,\CL,\CE)$ which 
is a generalization of the graph $C^*$-algebras.

The graph $C^*$-algebras $C^*(E)$ associated to (directed) graphs $E$ 
has drawn much interest  of many authors since 
they were introduced in \cite{KPRR} 
(\cite{FW} for finite graphs) and
has been generalized in many ways for the past twenty years. 
A graph $C^*$-algebra $C^*(E)$ is generated by a universal   
Cuntz-Krieger $E$-family consisting of 
projections $p_v$'s  and 
partial isometries $s_e$'s indexed by the vertices $v$ 
and the edges $e$ of $E$ satisfying 
a set of relations determined by the graph $E$ 
so that  every $C^*$-algebra generated by a 
Cuntz-Krieger $E$-family is a quotient of $C^*(E)$. 
By virtue of \cite{DT},  
it is enough to consider row-finite graphs 
as far as the properties  
that are preserved under strong Morita equivalence  are concerned. 
(IH) is one of such properties, and 
the following is known for a  row-finite graph $E$:
 
\begin{thm}\label{thm graph-ih} {\rm (\cite[Proposition 5.3]{BPRS} 
 and \cite[Theorem 3.9]{KPR})}
$C^*(E)$  has the property (IH),  
if and only if   $E$ satisfies condition (L)    
and every vertex connects to a loop.  
\end{thm}
 
\noindent 
Condition (L) means that every loop has an exit, and 
the theorem  actually holds true for arbitrary graphs 
(\cite[Corollary 2.14]{DT}). 
On the other hand, characterizations of 
purely infinite graph $C^*$-algebras 
were obtained in \cite{H, JP} 
for locally finite graphs with no sinks and 
in \cite{HZ} for arbitrary graphs. 
One consequence of \cite{HZ} that we are interested in 
is the following:

\begin{thm}\label{thm graph-pi} {\rm (\cite[Theorem 2.3]{HZ})} 
For a graph $E$, the following are equivalent: 
\begin{enumerate} 
\item[(a)] Every quotient of $C^*(E)$  has the property (IH).
\item[(b)] $C^*(E)=C^*(p_v,s_e)$ is purely infinite. 
\item[(c)]  There are no breaking vertices, and 
every projection $p_v$ is properly infinite.
\end{enumerate} 
\end{thm}

\noindent 
(If $E$ is row-finite, there are no breaking verices.) 
In the same paper, it is also shown that 
if $C^*(E)$ is purely infinite, then  
$E$ satisfies condition  (K) (that is, 
every vertex lies on no loops or lies on at least two loops  
each of which is not an initial path of the other) 
and that $E$ satisfies condition (K) 
if and only if 
every ideal of $C^*(E)$ is  invariant for the 
gauge action of $\mathbb T$ induced by universal property of 
$C^*(E)$.  

\vskip 1pc 

In order to obtain characterizations for 
the property (IH) or pure infiniteness of labeled 
graph $C^*$-algebras as in Theorem~\ref{thm graph-ih} 
or Theorem~\ref{thm graph-pi}, 
we need labeled space analogues to the notions used in the theorems. 
A graph with condition (L) was generalized to 
a {\it disagreeable} labeled space in \cite{BP2} and 
this notion will play an important role throughout this paper 
together with a notion of {\it loop} introduced in \cite{JKK}. 
Besides, we will heavily use the notion of {\it quotient labeled space} 
to deal with quotient algebras of labeled graph algebras.

We begin in Section 2 with the definitions of a labeled space and 
its $C^*$-algebra, and then review known facts about 
(quotient) labeled graph $C^*$-algebras and purely 
infinite $C^*$-algebras.
Briefly a labeled graph is a graph $E$ equipped  
with a labeling map $\CL:E^1\to \CA$ assigning a letter  
in the alphabet $\CA$ to each edge of $E$. 
If $\CB$ is an accommodating set, a certain collection of 
vertex subsets, 
then a triple $(E,\CL,\CB)$ is called a {\it labeled space}. 
It is known in \cite{BP1} that a $C^*$-algebra 
$C^*(E,\CL,\CB)=C^*(p_A, s_a)$, $A\in \CB$, $a\in \CA$, 
can be associated to $(E,\CL,\CB)$ in a similar way 
as a graph $C^*$-algebra could be to a graph.   
A labeled graph $C^*$-algebra $C^*(E,\CL,\CB)$, by definition,  
depends on the choice of an accommodating set $\CB$, and  
we will consider mostly the smallest accommodating set which 
we denote by $\CE$. 
Throughout this paper, 
we deal only  with labeled spaces $(E,\CL,\CE)$ 
of graphs $E$ with no sinks 
that are weakly left-resolving, set-finite, 
receiver set-finite, and normal unless stated otherwise.
The class of labeled graph $C^*$-algebras  includes and strictly includes 
all the graph $C^*$-algebras;  
for example, it is recently known \cite{JKKP} that 
there exists a family of unital simple 
labeled graph $C^*$-algebras 
which are A$\mathbb T$ algebras with nonzero $K_1$-groups, 
whereas any simple graph $C^*$-algebra is known to be 
either AF or purely infinite. 

\vskip 1pc 

In Section 3, we investigate the property (IH) 
for labeled graph $C^*$-algebras.   
Our main result  is as follows: 

\begin{thm} {\rm (Theorem~\ref{thm-IH}.(a))} 
If $(E,\CL,\CE)$ is a disagreeable labeled space such that 
every vertex  connects  to a loop, then $C^*(E,\CL,\CE)$ 
has the property (IH).  
\end{thm}
 
\noindent Note here that a loop in a labeled space is different from 
a cycle used in \cite{COP} while every cycle is a loop  and 
they are the same in a (trivially labeled) graph. 
In fact, we provide an example of 
a labeled space $(E,\CL,\CE_\om)$ 
(see Example~\ref{counterex}) without any cycles 
but with loops for which the $C^*$-algebra 
$C^*(E,\CL,\CE_\om)$ is simple and purely infninite. 
This example shows an obvious contrast between 
loops and cycles in a labeled space.
The converse of the theorem is also shown in 
Theorem~\ref{thm-IH}.(b) under some extra conditions.

\vskip 1pc 

In Section 4, we consider pure infiniteness of 
labeled graph $C^*$-algebras. 
We say that a labeled space $(E,\CL,\CE)$ 
is {\it strongly disagreeable} if
every quotient labeled space is disagreeable. 
Then we show in Lemma~\ref{GI} that 
if $(E,\CL,\CE)$ is strongly disagreeable, then 
every ideal of $C^*(E,\CL,\CE)$ is gauge-invariant. 
Thus the condition of being strongly disagreeable 
can be considered as an analogue for condition (K). 
The following is our main result 
for pure infiniteness of 
$C^*(E,\CL,\CE)$.

\begin{thm} {\rm (Theorem~\ref{SDPI})} 
Let $(E,\CL,\CE)$ be strongly disagreeable labeled space. 
Then the following are equivalent for $C^*(E,\CL,\CE)=C^*(s_a,p_A)$ : 
\begin{enumerate} 
\item[(a)]  Every quotient of $C^*(E,\CL,\CE)$ has the property (IH).
\item[(b)]  $C^*(E,\CL,\CE)$ is purely infinite.
\item[(c)] Every nonzero projection $p_A$ is properly infinite 
for $A\in \CE$.
\end{enumerate}
\end{thm}
 
\noindent 
We believe that  the condition of $(E,\CL,\CE)$ being strongly disagreeable 
should be a necessary condition for its $C^*$-algebra to be 
purely infinite. 
As an evidence, in Proposition~\ref{prop-nondisagreeable}, 
we show under some coditions 
that if $C^*(E,\CL,\CE)$ is purely infinite, then 
the labeled space $(E,\CL,\CE)$ is strongly disagreeable. 

\section{Preliminaries} 
\subsection{\bf Directed graphs and labeled spaces}

A  {\it directed graph} is a quadruple $E=(E^0,E^1,r,s)$
consisting of two countable sets of
vertices  $E^0$ and edges $E^1$,
and the range, source maps $r$, $s: E^1\to E^0$. 
If a vertex $v\in E^0$ emits (respectively, receives) no edges, 
it is called a {\it sink} (respectively, a {\it source}). 
A graph $E$ is {\it row-finite} if 
every vertex emits only finitely many edges.

For each $n\geq 1$, a finite path $\ld$ of {\it length} $n$, $|\ld|=n$, 
is a sequence of edges $\ld_{i}\in E^1$ such that 
$r(\ld_{i})=s(\ld_{i+1})$ for $ 1\leq i\leq n-1$. 
$E^n$  denotes the set of all  paths of length $n$, and 
the vertices in $E^0$ are regarded as finite paths of length zero. 
The maps $r,s$ naturally extend to the set 
$E^*=\cup_{n \geq 0} E^n$ of all finite paths, 
especially with $r(v)=s(v)=v$ for $v \in E^0$. 
We write $E^{\infty}$ for the set of all infinite paths 
$x=\ld_{1}\ld_{2}\cdots$, $\ld_{i}\in E^1$ with 
$r(\ld_i)=s(\ld_{i+1})$ for $i\geq 1$, 
and define $s(x):= s(\ld_1)$. 
For $A,B\subset E^0$ and $n\geq 0$, set
 $$ AE^n: =\{\ld\in E^n :  s(\ld)\in A\},\ \
  E^nB: =\{\ld\in E^n : r(\ld)\in B\},$$
 and  $AE^nB: =AE^n\cap E^nB$. 
 We write $E^n v$ for $E^n\{v\}$ and $vE^n$ for $\{v\}E^n$,
 and will use notation like $AE^{\geq k}$ and $vE^\infty$
 which should have their obvious meaning.

A finite path $\ld\in E^{\geq 1}$ with $r(\ld)=s(\ld)$ is 
called a {\it loop}, and  
an {\it exit} of a loop $\ld$ is a path 
$\dt\in E^{\geq 1}$ such that  
$|\dt|\leq |\ld|,\ s(\dt)=s(\ld), \text{ and } \dt\neq \ld_{[1,|\dt|]}.$ 
A graph $E$ is said to satisfy {\it condition} (L)  
if every loop has an exit and $E$ is said to satisfy {\it condition} (K) 
if no vertex in $E$ is the source vertex of exactly 
one loop which does not return to its source vertex more than once.

Let $\CA$ be a countable alphabet and let  
$\CA^*$ (respectively, $\CA^\infty$) denote   
the set of all finite words (respectively, infinite words)
in symbols of $\CA$. 
A {\it labeled graph} $(E,\CL)$ over $\CA$ consists of 
a directed graph $E$ and  a {\it labeling map} 
$\CL:E^1\to \CA$ which is always assumed to be onto. 
Given a graph $E$, one can define a so-called 
{\it trivial labeling} map 
 $\CL_{id}:=id:E^1\to E^1$ which is the identity map 
 on $E^1$ with the alphabet $E^1$. 
To each finite path $\ld=\ld_1\cdots \ld_n\in E^n$ 
of a labeled graph $(E,\CL)$ over $\CA$, there corresponds 
 a finite labeled path 
 $\CL(\ld):=\CL(\ld_1)\cdots \CL(\ld_n)\in \CL(E^n)\subset \CA^*$, and  
 similarly an infinite labeled path 
 $\CL(x):=\CL(\ld_1)\CL(\ld_2)\cdots\in \CL(E^\infty)\subset\CA^\infty$  
 to each infinite path $x=\ld_1\ld_2\cdots \in E^\infty$.  
 We often call these labeled paths just paths for convenience if 
 there is no risk of confusion,  
and use notation $\CL^*(E):=\CL(E^{\geq 1})$, 
where $E^{\geq 1}=E^*\setminus E^0$.
We also write $\CL(v):=v$ for $v\in E^0$ and 
$\CL(A):=A$ for $A\subset E^0$.  
 A subpath  $\af_i\cdots \af_j$ of  
$\af=\af_1\af_2\cdots\af_{|\af|}\in \CL^*(E)$ 
 is denoted by 
 $\af_{[i,j]}$ for $1\leq i\leq j\leq |\af|$, and 
 each $\af_{[1,j]}$, $1\leq j\leq |\af|$, is called an 
{\it initial path} of $\af$. 
 The range and source of a path $\af\in \CL^*(E)$ are defined by
 \begin{align*}
r(\af) &=\{r(\ld) \in E^0 \,:\, \ld\in E^{\geq 1},\,\CL(\ld)=\af\},\\
 s(\af) &=\{s(\ld) \in E^0 \,:\, \ld\in E^{\geq 1},\, \CL(\ld)=\af\},
\end{align*}
and the {\it relative range of $\af\in \CL^*(E)$  
with respect to $A\subset  E^0$} is defined by
$$
 r(A,\af)=\{r(\ld)\,:\, \ld\in AE^{\geq 1},\ \CL(\ld)=\af \}.
$$
A collection  $\CB$ of subsets of $E^0$ is said to be
 {\it closed under relative ranges} for $(E,\CL)$ if 
$r(A,\af)\in \CB$ whenever 
 $A\in \CB$ and $\af\in \CL^*(E)$. 
We call $\CB$ an {\it accommodating set}~ for $(E,\CL)$
 if it is closed under relative ranges,
 finite intersections and unions and contains 
$r(\af)$ for all $\af\in \CL^*(E)$.
A set $A\in \CB$ is called {\it minimal} (in $\CB$)  
if $A \cap B$ is either $A$ or $\emptyset$ for all $B \in \CB$.
By $\CB_{\rm min}$, we denote the set
$$ \CB_{\rm min}:=\{ A\in \CB:\, A\neq \emptyset\text{ and } 
A\cap B \text{ is either $A$ or }\emptyset 
\text{ for all  } B\in \CB\,\}$$ 
of all nonempty {\it minimal} sets in $\CB$.

If $\CB$ is accommodating for $(E,\CL)$, 
the triple $(E,\CL,\CB)$ is called
 a {\it labeled space}. 
We say that a labeled space $(E,\CL,\CB)$ is {\it set-finite}
 ({\it receiver set-finite}, respectively) if 
for every $A\in \CB$ and $k\geq 1$ 
 the set  $\CL(AE^k)$ ($\CL(E^k A)$, respectively) is finite.
  A labeled space $(E,\CL,\CB)$ is said to be {\it weakly left-resolving} 
if it satisfies 
 $$r(A,\af)\cap r(B,\af)=r(A\cap B,\af)$$
  for all $A,B\in \CB$ and  $\af\in \CL^*(E)$.
If $\CB$ is closed under relative complements, 
we call $(E,\CL, \CB)$ a {\it normal} labeled space as in \cite{BCP}. 

\vskip 1pc 

\noindent 
{\bf Assumption.}   Throughout this paper, we assume that 
graphs $E$ have no sinks and labeled spaces $(E,\CL,\CB)$ are  
weakly left-resolving,  set-finite, receiver set-finite, and 
normal unless stated otherwise. 

\vskip 1pc

By $\CE$  we denote the smallest accommodating set 
for which $(E,\CL, \CE)$ is a normal labeled space. 
 Let  $\Om_0(E)$ be the set of all vertices that are not sources, 
and for each $l\geq 1$, define 
a relation $\sim_l$ on  $\Om_0(E)$ by 
$v\sim_l w$ if and only of $\CL(E^{\leq l} v)=\CL(E^{\leq l} w)$.
Then $\sim_l$ is an equivalence relation, and the equivalence class 
$[v]_l$ of $v\in \Om_0(E)$ is called a {\it generalized vertex} 
(or simply a vertex if there is no risk of confusion).  
   If $k>l$,  then $[v]_k\subset [v]_l$ is obvious and
   $[v]_l=\cup_{i=1}^m [v_i]_{l+1}$
   for some vertices  $v_1, \dots, v_m\in [v]_l$ 
(\cite[Proposition 2.4]{BP2}). 
The generalized vertices of labeled graphs 
play the role of vertices in usual graphs.   
Moreover,  we have
\begin{eqnarray}\label{CE}
\CE=\Big\{ \cup_{i=1}^n [v_i]_l:\, v_i\in E^0,\  l\geq 1,\,  n\geq 0 \Big\},
\end{eqnarray}
with the convention  $\sum_{i=1}^0 [v_i]_l:=\emptyset$ 
by \cite[Proposition 2.3]{JKK}. 

\vskip 1pc 

\subsection{$C^*$-algebras of labeled spaces}

We review the definition of  $C^*$-algebras 
associated to labeled spaces from \cite{BP1, BP2}.

\begin{dfn} 
\label{def-representation}
A {\it representation} of a labeled space $(E,\CL,\CB)$
is a family of projections $\{p_A\,:\, A\in \CB\}$ and
partial isometries
$\{s_a\,:\, a\in \CA\}$ such that for $A, B\in \CB$ and $a, b\in \CA$,
\begin{enumerate}
\item[(i)]  $p_{\emptyset}=0$, $p_{A\cap B}=p_Ap_B$, and
$p_{A\cup B}=p_A+p_B-p_{A\cap B}$,
\item[(ii)] $p_A s_a=s_a p_{r(A,a)}$,
\item[(iii)] $s_a^*s_a=p_{r(a)}$ and $s_a^* s_b=0$ unless $a=b$,
\item[(iv)]\label{CK4}  $p_A=\sum_{a\in \CL(AE^1)} s_a p_{r(A,a)}s_a^*.$ 
\end{enumerate}
\end{dfn}

\vskip 1pc

\noindent
It is known \cite[Theorem 4.5]{BP1} that 
given a labeled space  $(E,\CL,\CB)$, 
there exists a $C^*$-algebra $C^*(E,\CL,\CB)$ generated by 
a universal representation $\{s_a,p_A\}$ of $(E,\CL,\CB)$. 
We call $C^*(E,\CL,\CB)$ the {\it labeled graph $C^*$-algebra} of
a labeled space $(E,\CL,\CB)$ which is unique up to isomorphism, 
and simply write $C^*(E,\CL,\CB)=C^*(s_a,p_A)$ 
to indicate the generators. 
Note that   $s_a\neq 0$ and $p_A\neq 0$ for $a\in \CA$
and  $A\in \CB$, $A\neq \emptyset$. 
Every graph $C^*$-algebra    
is a labeled graph $C^*$-algebra of a trivial labeled space. 
The labeled graph $C^*$-algebra $C^*(E,\CL,\CB)$ 
depends on the choice of an accommodating set $\CB$, but 
we are mainly interested in $C^*$-algebras of the labeled spaces 
$(E,\CL,\CE)$ throughout this paper.

\vskip 1pc

\begin{remark}\label{basics} 
Let $(E,\CL,\CB)$ be a labeled space with $C^*(E,\CL,\CB)=C^*(s_a,p_A)$. 
By $\ep$, we denote a symbol (not belonging to $\CL^*(E)$)
such that $r(\ep)=E^0$ and $r(A,\ep)=A$ for all $A \subset E^0$. 
We write $\CL^{\#}(E)$ for the union $\CL^*(E)\cup \{\ep\}$.  
Let $s_\ep$ denote the unit of the multiplier algebra of $C^*(E,\CL,\CB)$.  

\begin{enumerate}
\item[(a)] 
We have the following equality  
$$(s_\af p_{A} s_\bt^*)(s_\gm p_{B} s_\dt^*)=
\left\{
   \begin{array}{ll}
      s_{\af\gm'}p_{r(A,\gm')\cap B} s_\dt^*, & \hbox{if\ } \gm=\bt\gm' \\
      s_{\af}p_{A\cap r(B,\bt')} s_{\dt\bt'}^*, & \hbox{if\ } \bt=\gm\bt'\\
      s_\af p_{A\cap B}s_\dt^*, & \hbox{if\ } \bt=\gm\\
      0, & \hbox{otherwise,}
   \end{array}
\right.
$$                    
for $\af,\bt,\gm,\dt\in  \CL^{\#}(E)$ and $A,B\in \CB$ 
(see \cite[Lemma 4.4]{BP1}). 
Since 
$s_\af p_A s_\bt^*\neq 0$ if and only if 
$A\cap r(\af)\cap r(\bt)\neq \emptyset$, 
it follows that  
\begin{eqnarray}\label{eqn-elements}
\hskip 3pc 
C^*(E,\CL,\CB)=\overline{\rm span}\{s_\af p_A s_\bt^*\,:\,
\af,\,\bt\in  \CL^{\#}(E) ~\text{and}~ A \subseteq r(\af)\cap r(\bt)\}. 
\end{eqnarray}  

\item[(b)] Universal property of  $C^*(E,\CL,\CB)=C^*(s_a, p_A)$ 
defines a strongly continuous action 
$\gm:\mathbb{T} \rightarrow Aut(C^*(E,\CL,\CB))$,  
called the  {\it gauge action}, such that
$$\gm_z(s_a)=zs_a ~ \text{ and } ~ \gm_z(p_A)=p_A$$
for $a\in \CA$ and $A\in \CB$. 
Averaging over $\gm$ with respect to the normalized Haar measure 
of the compact group  $\mathbb{T}$, 
$$\Phi(a):=\int_{\mathbb{T}}\gm_z(a)dz,\ a \in C^*(E,\CL,\CB),$$
defines 
a conditional expectation $\Phi:C^*(E,\CL,\CB)\to C^*(E,\CL,\CB)^{\gm}$
onto the fixed point algebra which is 
also known to be faithful.

\item[(c)] It is shown in \cite[Theorem 4.4]{BP2} that 
$C^*(E,\CL,\CB)^{\gm}$ is an AF algebra 
isomorphic to 
$C^*(E,\CL,\CB)^{\gm} \cong \overline{\cup_{k,l}(\oplus_{[v]_l}\CF^k([v]_l)})$, 
where 
\begin{eqnarray}\label{Fk}
\CF^k([v]_l)=\overline{\rm span}\{s_{\af}p_{[v]_l}s_{\bt}^*:\af,\bt \in \CL(E^k)\}.
\end{eqnarray}

\end{enumerate}
\end{remark}

\vskip 1pc
\subsection{Cuntz-Krieger uniqueness theorem}
Recall that a {\it Cuntz-Krieger $E$-family} for a graph $E$ 
is a representation of the labeled space 
$(E,\CL_{id},\CE)$ with the trivial labeling, and 
{\it the Cuntz-Krieger uniqueness theorem} for graph $C^*$-algebras 
says that 
if $E$ satisfies condition (L), 
then every Cuntz-Krieger $E$-family of nonzero 
operators generates the same $C^*$-algebra $C^*(E)$ 
up to isomorphism  
(for example, see \cite[Theorem 3.1]{BPRS}, \cite[Corollary 2.12]{DT}, 
and  \cite[Theorem 3.7]{KPR}).  
A condition for  labeled spaces corresponding 
to condition (L) for directed graphs was given in \cite[Definition 5.2]{BP2}, 
and  we briefly review  it here. 

A labeled path $\af \in \CL^*(E)$ with 
$s(\af)\cap [v]_l \neq \emptyset$ is 
called {\it agreeable} for $[v]_l$ if 
$\af=\bt\af'=\af'\gm$ for some $\af',\bt,\gm \in \CL^*(E)$ 
with $|\bt|=|\gm| \leq l$. 
Otherwise $\af$ is called {\it disagreeable}.
(Note that any  path  $\af$ agreeable for $[v]_l$ 
must be of the form $\af=\bt^k\bt'$ for 
some $\bt\in E^{\leq l}$,  $k\geq 0$, 
and an initial path $\bt'$ of $\bt$.) 
 We say that $[v]_l$ is {\it disagreeable} 
 if there is an $N\geq 1$ 
 such that for all $n > N$ there is an 
 $\af \in \CL(E^{ \geq n})$ which is disagreeable for $[v]_l$. 
 (It can be easily seen that 
 $[v]_l$ is not disagreeable if and only if 
 there is an $N\geq 1$ such that 
 every path in $\CL([v]_l E^{ \geq N})$ is  agreeable for $[v]_l$.)
A labeled space $(E,\CL, \CB)$ is said to be  {\it disagreeable}  
 if for every $v\in E^0$, there is an $L_v\geq 1$ such that 
  every $[v]_l$ is disagreeable for all $l \geq L_v$. 
Then it is shown in \cite[Lemma 5.3]{BP2} that 
for a graph $E$, 
the labeled space $(E,\CL_{id}, \CE)$ with the trivial labeling 
is disagreeable if and only if $E$ satisfies the condition (L). 

In \cite[Definition 9.5]{COP}, 
the notion of cycle was introduced to define 
condition $(L_\CB)$ for a labeled space $(E,\CL,\CB)$ 
(more generally for 
Boolean dynamical systems) which can be regarded as 
another condition equivalent to condition  (L) for usual directed graphs. 
For a path $\af\in \CL^*(E)$ and a  set $\emptyset\neq A\in \CE$, 
the pair $(\af, A)$ is called a {\it cycle} if 
$B=r(B,\af^k)$ holds for all $k\geq 0$ 
and  nonempty subsets  $B\in \CE$ of $A$ 
(see the proof of \cite[Proposition 9.6]{COP}). 
We say that a labeled space $(E,\CL,\CE)$ satisfies  
condition $(L_\CE)$ if there is no cycle without exits.
We will see in Proposition~\ref{QLE} that 
every disagreeable labeled space satisfies condition $(L_\CE)$. 
The following is the {\it Cuntz-Krieger uniqueness theorem} 
for labeled graph $C^*$-algebras.
\vskip 1pc 

\begin{thm} {\rm (\cite[Theorem 5.5]{BP1}, \cite[Theorem 9.9]{COP})} 
\label{CK uniqueness thm} 
Let  $\{t_a, q_A\}$ be a representation of a labeled space $(E,\CL,\CE)$ 
such that $q_A\neq 0$ for all nonempty $A\in \CE$. 
If $(E,\CL,\CE)$ satisfies condition  $(L_\CE)$, in particular 
if  $(E,\CL,\CE)$ is disagreeable, then 
the canonical homomorphism $\phi:C^*(E,\CL,\CE)=C^*(s_a, p_A)\to 
 C^*(t_a, q_A)$ such that $\phi(s_a)=t_a$ and 
$\phi(p_A)=q_A$  is an isomorphism. 
\end{thm} 

\vskip 1pc
   
\subsection{Ideal structure of labeled graph $C^*$-algebras} 
We first review definitions of quotient labeled spaces and their $C^*$-algebras 
which were introduced in \cite{JKP} to 
study the ideal structure of labeled graph $C^*$-algebras. 
 
Let $(E,\CL,\CE)$ be a labeled space and 
$\sim_R$ an equivalence relation on $\CE$. 
Denote the equivalence class of $A \in \CE$ by $[A]$ 
(rather than $[A]_R$) and set  
\begin{align*}
\CA_R:& =\{a \in \CA :[r(a)] \neq [\emptyset]\}\\ 
\CL_R(E^n):& 
   =\{\af\in \CL(E^n): [r(\af)]\neq [\emptyset]\} \text{ for }\ n\geq 1,\\ 
\CL_R^*(E): & =\{\af \in \CL^*(E): [r(\af)] \neq [\emptyset]\}\\
\CL_R([A]E^n):& 
   =\{\af \in \CL([A]E^n): [r(\af)] \neq [\emptyset]\} \text{ for }\ n\geq 1,
\end{align*}
and so forth.
If the following operations
$$[A]\cup[B]:=[A\cup B], \ \ [A]\cap[B]:=[A\cap B],
\ \  [A]\setminus[B]:=[A\setminus B]$$
are well-defined on $\CE/R:=\{ [A]: A\in \CE\}$ and if 
the relative range 
$$r([A],\af):=[r(A,\af)]$$ 
is well-defined for $[A]\in \CE/R$ and $\af \in \CL_R^*(E)$ 
so that 
$r([A],\af)=[\emptyset]$ for all $\af \in \CL_R^*(E)$ implies 
$[A] = [\emptyset]$, 
we call $(E,\CL,\CE/R)$ a {\it quotient labeled space} of 
$(E,\CL,\CE )$. 

 We say that a quotient labeled space $(E,\CL,\CE/R)$ is 
 {\it weakly left-resolving} if 
 for  $[A],[B] \in \CE/R$ and $\af \in \CL_R^*(E)$, 
 the following holds:
 $$r([A],\af) \cap r([B], \af)=r([A]\cap[B],\af).$$  
We write  $[A] \subseteq [B]$ if $[A] \cap [B]=[A]$, and 
$[A] \subsetneq [B]$ if $[A] \subseteq[B]$ and  $[A] \neq [B]$.

\vskip 1pc

\begin{dfn}\label{QR} (\cite[Definition 3.3]{JKP}) 
 A {\it representation}  of a quotient labeled space 
 $(E,\CL,\CE/R)$ (always assumed weakly left-resolving) is a family of 
 projections $\{p_{[A]}: [A] \in \CE/R\}$ and partial isometries $\{s_a: a \in \CA_R\}$ 
 such that
\begin{enumerate}
\item[(i)]  $p_{[\emptyset]}=0$, $p_{[A]\cap [B]}=p_{[A]}p_{[B]}$, and
$p_{[A]\cup [B]}=p_{[A]}+p_{[B]}-p_{[A]\cap [B]}$,
\item[(ii)] $p_{[A]} s_a=s_a p_{r([A],a)}$,
\item[(iii)] $s_a^*s_a=p_{[r(a)]}$ and $s_a^* s_b=0$ unless $a=b$,
\item[(iv)] $ p_{[A]}=\sum_{a\in \CL_R([A]E^1)} s_a p_{r([A],a)}s_a^*\ $ 
 if $\CL_R([A]E^1) \neq \emptyset$. 
\end{enumerate}
\end{dfn}

 \vskip 1pc
  
\noindent 
It is known \cite[Theorem 3.10]{JKP} that if $(E, \CL, \CE/R)$ is 
a  quotient labeled space,  there exists a $C^*$-algebra  
$C^*(E,\CL,\CE/R)$, called the {\it quotient labeled graph $C^*$-algebra} of 
$(E,\CL,\CE/R)$, 
generated by a universal representation of $(E,\CL,\CE/R)$. 
By \cite[Lemma 3.8]{JKP}, we have 
\begin{eqnarray}\label{q-span}
C^*(E,\CL,\CE/R)=\overline{\rm span}\big\{ s_\af p_{[A]} s_\bt^*: 
\af,\bt\in \CL_R^*(E),\ [A]\neq [\emptyset]\,\big\}.
\end{eqnarray}

The equivalence relation on $\CE$ which we are most interested in is 
related with  hereditary saturated subsets $H$ of $\CE$.  
Recall  from \cite[Definition 3.4]{JKP} 
that a subset  $H$ of $\CE$  is {\it hereditary} if 
it is closed under finite unions, relative ranges, 
and subsets in $\CE$.
A hereditary set $H$ is {\it saturated} if  $A \in H$ 
whenever $A \in \CE $ satisfies $r(A,\af) \in H$ 
for all $\af \in \CL^*(E)$. 
Every hereditary saturated subset $H\subset \CE$ defines 
an equivalence relation  $\sim_H$ on $\CE$; 
for $A, B\in \CE$,
$$A\sim_H B \Longleftrightarrow A\cup W=B\cup W\ \text{ for some } W\in H.$$ 
For a hereditary set $H$, let  
$\overline{H}$ be the smallest 
hereditary saturated set containing $H$. 
Then the ideal $I_H$ of $C^*(E, \CL,\CE )$ generated by the projections 
$\{p_A: A \in H\}$ is gauge-invariant (\cite[Lemma 3.7]{JKP}) 
and 
$$I_H 
=\overline{\rm span}
 \{s_{\mu}p_As_{\nu}^*: \mu,\nu \in \CL^\#(E),\, A \in \overline{H} \}.
$$

\vskip 1pc
 
\begin{remark} \label{remark-review} 
Let $(E,\CL,\CE)$ be a labeled space.
\begin{enumerate}
\item[(a)]  $H\mapsto I_H$ is an inclusion preserving 
bijection between the nonempty hereditary saturated subsets of 
$\CE$ and the nonzero gauge-invariant ideals of 
the $C^*$-algebra $C^*(E,\CL,\CE)$. 
Moreover, the quotient algebra 
$C^*(E,\CL,\CE)/I_H$  is isomorphic to 
$C^*(E,\CL, \CE/H)$ (\cite[Theorem 5.2]{JKP}). 

\item[(b)]  
 Every quotient labeled graph $C^*$-algebra 
$C^*(E, \CL,\CE/R)=C^*(s_a,p_{[A]})$ admits the {\it gauge action} $\gm$ 
of $\mathbb{T}$, 
$$\gm_z(s_a)=zs_a ~\text{and}~ \gm_z(p_{[A]})=p_{[A]}$$
for $a \in \CA_R$ and $[A] \in \CE/R$. 
 It can be seen as in \cite[Theorem 4.4]{BP2} that 
 the fixed point algebra  
$$
\hskip 3pc C^*(E,\CL,\CE/R)^{\gm}=
\overline{\rm span}\{s_{\af}p_{[A]}s_{\bt}^* : \af, \bt  \in \CL_R(E^n), 
   \, n\geq 1,\,  [A]\in \CE/R \} 
$$
is  an  AF algebra (see the proof of \cite[Theorem 4.2]{JKP}). 
Also, one sees as in Remark~\ref{basics}.(b) 
that there exists a faithful conditional expectation 
$\Phi$ of $C^*(E,\CL,\CE/R)$ onto the fixed point algebra 
$C^*(E,\CL,\CE/R)^{\gm}$.
\end{enumerate}
\end{remark}

\vskip 1pc

\subsection{Simplicity of labeled graph $C^*$-algebras}

Let $(E,\CL,\CE)$ be a labeled space and let
$$\overline{\CL(E^\infty)}:=\{x\in \CA^{\mathbb N}\mid 
x_{[1,n]}\in \CL(E^n) \ \text{for all } n\geq 1\}$$ 
be the set of all  infinite sequences $x$ such that 
every finite subpath of $x$ occurs as a labeled path in $(E,\CL)$. 
Clearly $\CL(E^\infty)\subset \overline{\CL(E^\infty)}$, and 
in fact, $\overline{\CL(E^\infty)}$ is the closure of $\CL(E^\infty)$ in the 
totally disconnected perfect space $\CA^{\mathbb N}$ 
which has the topology with  a countable 
basis of open-closed cylinder sets 
$Z(\af):=\{x\in \CA^{\mathbb N}: x_{[1,n]}=\af\}$, $\af\in \CA^n$, $n\geq 1$ 
(see Section 7.2 of \cite{Ki}). 
 We say that a labeled space $(E,\CL,\CE )$ is 
{\it strongly cofinal} if for each $x\in \overline{\CL(E^\infty)}$ 
and $[v]_l\in \CE $, there exist 
an $N\geq 1$ and a finite number of paths 
$\ld_1, \dots, \ld_m\in \CL^*(E)$ such that 
$$r(x_{[1,N]})\subset \cup_{i=1}^m r([v]_l,\ld_i).$$ 

\vskip 1pc

\begin{remarks}\label{simple condition} 
For simplicity of $C^*(E,\CL,\CE)$, we 
note the following. 
\begin{enumerate}
\item[(a)] If $(E, \CL, \CE )$ is disagreeable and strongly cofinal, 
then $C^*(E,\CL,\CE )$ is simple 
(see \cite[Theorem 3.16]{JK} and \cite[Remark 3.8]{JKKP}):  
a condition of being  disagreeable  and cofinal
assumed in \cite[Theorem 6.4]{BP2} is in fact not sufficient 
for the $C^*$-algebra $C^*(E,\CL,\CE)$ to be simple, 
but modifying its proof for disagreeable and strongly cofinal 
labeled space one can obtain the simplicity result.
 
\item[(b)] It is recently obtained in \cite[Theorem 9.16]{COP}
that  $C^*(E,\CL,\CE )$ is simple if and only if 
$(E,\CL,\CE)$ satisfies condition $L_\CE$ and 
the only hereditary and saturated subsets of $\CE$  are 
$\emptyset$ and $\CE$. 
\end{enumerate}

\end{remarks}

\vskip 1pc
 
\subsection{Purely infinite $C^*$-algebras}  
For positive elements $a, b$ in a $C^*$-algebra $A$, 
we write $a \precsim b$  if 
there exists a sequence $\{x_k\}_{k=1}^{\infty}$ in $A$ 
such that $x_k^*bx_k \rightarrow a$. 
More generally, for positive elements $a \in M_n(A)$ and $b \in M_m(A)$, 
write $a \precsim b$  if 
there exists a sequence $\{x_k\}_{k=1}^{\infty}$ in $M_{m,n}(A)$ such that 
$x_k^*bx_k \rightarrow a$.
A positive element $a\in A$ is said to be {\it  infinite} if 
there exists a nonzero positive element $b \in A$ such that $a \oplus b \precsim a$, 
and {\it properly infinite} if $a \oplus a \precsim a$. 
Every properly infinite $a$ is obviously infinite. 
Because the relation $a \oplus a \precsim a$ is preserved under $*$-homomorphisms, 
if $a\in A$ is properly infinite, then for any ideal $I$ of $A$ 
the element $a+I$ is either zero or properly infinite. 

 A $C^*$-algebra $A$ is said to be 
 {\it purely infinite} if there are no characters on $A$ and 
 if for every pair of positive elements $a,b \in A$,  
  $a \precsim b$ if and only if $a \in \overline{AbA}$ 
 (see \cite[Definition 4.1]{KR}\label{PI}). 
 Note that any $C^*$-algebra of type $I$  is 
 not purely infinite (\cite[Definition 4.4]{KR}). 
 For example, $C(\mathbb T)\otimes M_n$ is not purely infinite for $n\geq 1$.

\begin{remarks}\label{rmk-pre}
Let $A$ be a $C^*$-algebra and $a\in A$. 
\begin{enumerate}
\item[(a)] $A$ is purely infinite if and only if 
every nonzero positive element in A is 
properly infinite (\cite[Theorem 4.16]{KR}). 
For projections  $p,q$ in A, we write $p \preceq q $ 
if $p$ is Murray-von Neumann equivalent to a subprojection of $q$, 
namely $p\sim q\leq p$.
A projection $p\in A$ is  infinite if and only if  $p\sim q\lneq p$ for 
some subprojection $q$ of $p$,
and is properly infinite if and only 
$p$ has mutually orthogonal subprojections $p_1$ and $p_2$ 
such that $p_1\sim p_2\sim p$. 

\item[(b)]  $A$ is purely infinite if 
every nonzero hereditary $C^*$-subalgebra in every quotient 
algebra of $A$ 
contains an infinite projection (\cite[Proposition 4.7]{KR}). 

\item[(c)] Every nonzero hereditary $C^*$-subalgebra and 
every quotient algebra of a purely infinite $C^*$-algebra 
is again purely infinite (\cite[Proposition 4.17 and Proposition 4.3]{KR}).

\item[(d)] If $A$ is simple, then $A$ is purely infinite 
if and only if it has the property 
(IH) (see \cite[Proposition 4.6 and Proposition 5.4]{KR}). 
It should be noted that for non-simple $C^*$-algebras, 
one of these two properties is 
neither weaker nor stronger than the other, in general 
(for example, see \cite[Example 4.6]{KR}).
\end{enumerate}
\end{remarks}
 
\vskip 1pc
 
\section{Labeled graph $C^*$-algebras with the property (IH)}

\subsection{Infinite $C^*$-algebras of disagreeable labeled spaces}
Because  the original definition of a disagreeable labeled space 
seems a bit complicated,  
we list its equivalent but simpler conditions in 
Proposition~\ref{prop-disagreeable} 
after we observe the following lemma which will be frequently used. 

\vskip 1pc 

\begin{lem}\label{path comparison}
Let $\af$ and $\bt$ be paths in $\CL^*(E)$ 
$($or words in $\CA^*$$)$ with $|\af|\leq |\bt|$. 
If $\af^m=\bt^n$ for some $m, n\geq 1$,
then $\af,\bt\in \{\af_0^k: k\geq 1\}$ 
for an initial path $\af_0$ of $\af$.
\end{lem} 
\begin{proof}
If $|\af|=1$ or $n=1$, there is nothing to prove.
Thus  we assume $|\af|>1$ and $n>1$.  
Since $|\af|\leq |\bt|$ and the path $\af^m=\bt^n$ begins  with 
both of $\af$ and $\bt$, the longer path $\bt$ must have $\af$ as its 
initial path. 
If $\bt=\af^{n_1}$ for some $n_1\geq 1$, we are done. 
Otherwise, we can write $\bt=\af^{n_1}\af^{(1)}$ for some 
$n_1\geq 1$ and an initial path 
$\af^{(1)}$ of $\af$ with $1\leq |\af^{(1)}|<|\af|$. 
Then by canceling $\af^{n_1}$  
from the left of both sides of $\af^m=\bt^n$, we obtain 
$\af^{m-n_1}=\af^{(1)}\bt^{n-1}$.  
Thus, with a subpath $\af^{(2)}$ of $\af$ such that 
$\af=\af^{(1)}\af^{(2)}$, 
it follows that  
$$\af^{(1)}\af^{(2)}\cdots \af^{(1)}\af^{(2)} 
=\af^{(1)}\af^{n_1}\af^{(1)}\cdots \af^{n_1}\af^{(1)}$$ 
in which the right hand side is equal to 
$\af^{(1)}(\af^{(1)}\af^{(2)})^{n_1}\af^{(1)}\cdots \af^{n_1}\af^{(1)}$, hence 
comparing initial parts of both sides we see that 
$\af^{(2)}\af^{(1)}=\af^{(1)}\af^{(2)}$. 
Thus the lemma reduces to the following claim: 

\vskip .5pc
\noindent
{\bf Claim}. 
If $\af,\bt\in \CL^*(E)$ and $\af\bt=\bt\af$, then 
$\af,\bt\in \{\dt^k: k\geq 1\}$ for some $\dt\in \CL^*(E)$. 

\vskip .5pc
\noindent 
Again we may assume that $|\af|<|\bt|$ and   
$\bt=\af^m \af^{(1)}$ for some $m\geq 1$ and an initial path 
$\af^{(1)}$ of $\af$ with $1\leq |\af^{(1)}|<|\af|$. 
Then $\af\bt=\bt\af$ gives $\af\af^{(1)}=\af^{(1)}\af$.
If $\af$ is a repetition of $\af^{(1)}$, that is, 
$\af=(\af^{(1)})^n$ for some $n\geq 1$, the claim follows. 
If not, with $\af=\af^{(1)}\af^{(2)}$, 
we have from  $\af\af^{(1)}=\af^{(1)}\af$  
that $\af^{(1)}\af^{(2)}=\af^{(2)}\af^{(1)}$. 
If one of $\af^{(1)}$ or $\af^{(2)}$ is a repetition of the other, 
for example if $\af^{(1)}= (\af^{(2)})^s$ for some $s\geq 1$, 
we can take $\af^{(2)}$ for $\dt$. 
If not, we continue the process to obtain 
\begin{eqnarray}\label{alpha_k} 
\af^{(k)}\af^{(k+1)}=\af^{(k+1)}\af^{(k)}
\end{eqnarray} 
for initial paths $\af^{(k)}$, $k\geq 1$, of $\af$. 
We can check if one of $\af^{(k)}$ or $\af^{(k+1)}$ is a repetition  of the other 
whenever we arrives at step (\ref{alpha_k}). 
If it is the case, then we stop and obtain the claim. 
If not, again we can go on to the next step but with paths having 
length smaller than the previous step. 
Obviously this process should stop after finite steps, and 
when it stops we obtain an initial path $\dt$ of $\af$ satisfying the claim.
\end{proof}

\vskip 1pc

\begin{prop} \label{prop-disagreeable}
For a labeled space $(E,\CL,\CB)$, the following are equivalent:
\begin{enumerate}
\item[(a)] $(E,\CL,\CB)$ is disagreeable.  
\item[(b)] $[v]_l$ is disagreeable for all $v\in E^0$ and $l\geq 1$.
\item[(c)] For each nonempty $A\in \CB$ and a path $\bt\in \CL^*(E)$, 
there is an $n\geq 1$ such that $\CL(AE^{|\bt|n})\neq \{\bt^n\}$.
\end{enumerate}
\end{prop}

\begin{proof} 
The first two are known to be equivalent 
in \cite[Proposition 3.9.(iii)]{JK}.

To show (b) $\Rightarrow$ (c), 
suppose  there are a nonempty set $A\in \CE$  and 
a path $\bt\in \CL^*(E)$ such that   
$\CL(AE^{|\bt|n})=\{\bt^n\}$ for all $n\geq 1$.  
Then for a vertex   $\emptyset\neq [v]_{l}\subset A$ with $l\geq |\bt|$,  
$\CL([v]_{l}E^{|\bt|n})=\{\bt^n\}$ for all $n\geq 1$, 
and $[v]_{l}$ is not disagreeable.

For (c) $\Rightarrow$ (b), 
suppose $[v]_l$ is not disagreeable for some $v\in E^0$ and $l\geq 1$.  
Then, by \cite[Proposition 3.9.(i)]{JK}, 
there is an $N\geq l$ such that 
every path in $\CL([v]_lE^{\geq N})$ is agreeable for $[v]_l$.  
Observe here that for any $\dt\in \CL^*(E)$, 
every path in $\CL(r([v]_l,\dt^N)E^{\geq 1})$ is agreeable 
for $[v]_l$,  
which is immediate from the inclusion relation 
$\CL(r([v]_l,\dt^N)E^{\geq 1})\subset \CL([v]_lE^{\geq N})$.
One can choose a (finite) set 
$\{\dt_1, \dots, \dt_m\}$ of paths 
$\dt\in \CL^*(E)$ with $|\dt|\leq l$ such that 
every path in $\CL([v]_l E^{\geq N})$ 
is of the form 
\begin{eqnarray}\label{form}
\dt_i^{n}\dt_i'
\end{eqnarray} for some $n\geq 1$, $1\leq i\leq m$,  and 
$\dt_i'$ an initial path  of $\dt_i$.
We may assume that 
\begin{eqnarray}\label{delta}
\text{if } i\neq j, \text{ then }|\dt_i^{n_i}|=|\dt_j^{n_j}|  
\text{ and }\dt_i^{n_i}\neq \dt_j^{n_j} 
\text{ for some }n_i, n_j\geq 1 .
\end{eqnarray}
In fact, if $\dt_1^{n_1}=\dt_2^{n_2}$ 
whenever $|\dt_1^{n_1}|=|\dt_2^{n_2}|$, 
then by Lemma~\ref{path comparison} 
we can find an initial path $\dt$ of $\dt_1$ 
such that  $\dt_1$ and $\dt_2$ are equal to 
some repetitions of $\dt$. 
Then we remove $\dt_1$, $\dt_2$ from the set 
$\{\dt_1, \dots, \dt_m\}$, and include $\dt$ instead.  
Pick $L\geq N$  large enough so that 
the path $\dt_1^L$ can not be of the form in 
(\ref{form}) for $i\neq 1$ 
(this is possible by (\ref{delta})).
Then for $A:=r([v]_l, \dt_1^L)\in \CE$, 
every path in $\CL(AE^{\geq 1})$, 
already observed as agreeable for $[v]_l$, 
must be of the form $\dt_1^n\dt_1'$. 
But this contradicts to (c).
\end{proof}

\vskip 1pc 

We can naturally extend the notion of `disagreeable  labeled space' to 
`disagreeable quotient labeled space' as follows. 
The Cuntz-Krieger uniqueness theorem for 
disagreeable quotient labeled spaces will be given
in Section 4. 

\vskip 1pc  
    
\begin{dfn}  \label{disagreeable q}
Let  $(E,\CL,\CE/R)$ be a quotient labeled space 
and $[[v]_l] \neq [\emptyset]$ in $\CE/R$. 
A path $\af \in \CL_R^*(E)$ with  
$s(\af)\cap [[v]_l]\neq \emptyset$
is said to be {\it agreeable} for $[[v]_l]$ if 
$$[r([v]_l,\af)]= [\emptyset]\ \text{ or } \  
\af = \bt\af'=\af'\gm$$ 
for some $\af', \bt,\gm \in \CL_R^*(E)$ 
with  $ |\bt|=|\gm| \leq l $. 
Otherwise $\af$ is  {\it disagreeable} for $[[v]_l]$.
Also  $[[v]_l]$ is called {\it disagreeable} if 
there is an $N \geq 1$ such that for all $ n > N$, 
$$\CL_R([[v]_l]E^{\geq n}):=\CL([[v]_l]E^{\geq n}) \cap \CA_R^*$$ 
has a disagreeable path for $[[v]_l]$. 
The quotient labeled space $(E,\CL,\CE/R)$ is said to be 
{\it disagreeable} if
every $[[v]_l] \neq [\emptyset]$ is disagreeable 
for $l\geq 1$ and $v \in E^0$.
\end{dfn}

\vskip 1pc  

\begin{remark}\label{q-disagreeable} 
As in Proposition~\ref{prop-disagreeable}, 
one can  prove that 
$(E,\CL,\CE/R)$ is disagreeable if and only if 
for each  $[A]\neq [\emptyset]$  in $\CE/R$ and 
a path $\bt\in \CL_R^*(E)$, 
there is an $n\geq 1$ such that $\CL_R([A]E^{|\bt|n})\neq \{\bt^n\}$.
\end{remark}
 
\vskip 1pc 

It is easy to see that if a graph $E$ has a loop with an exit, 
the $C^*$-algebra 
$C^*(E)$ contains an infinite projection. 
In order to see whether this is true  for  
labeled graph $C^*$-algebras,  
we will use  the following definitions of loop and exit 
in a quotient labeled space (see \cite[Definition 3.2]{JKK}).

\vskip 1pc 

\begin{dfn}\label{qloop}
Let $(E, \CL,  \CE/R)$ be a quotient labeled space. 
For a  path $\af \in \CL_R^*(E)$ and  
a set $\emptyset\neq [A] \in  \CE/R$, 
we call  $(\af,[A])$ a {\it loop} if  $[A] \subseteq r([A],\af)$.  
We say that a loop  $(\af,[A])$ has an {\it exit} 
if one of the following holds:
\begin{enumerate}
\item[(i)] there exists a path 
$\bt \in \CL_R([A]E^{\geq 1})$ such that 
$$|\bt|=|\af|,\ \bt \neq \af,\ \text{and } 
r([A],\bt) \neq [\emptyset],$$
\item[(ii)] $[A]\subsetneq r([A],\af)$.
\end{enumerate} 
Clearly every cycle is a loop, and 
if $(\af,A)$ is a cycle with an exit, the exit must be 
of type (i). 

\end{dfn}
 
\vskip 1pc 

\begin{ex}\label{ex-1} 
We give examples of labeled spaces with loops 
which have exits 
of types (i) and (ii), respectively.  
\begin{enumerate}
\item[(i)]
 For the  labeled graph $(E, \CL)$ given below, 
 let $H$ be the smallest hereditary saturated subset of $\CE$ 
 containing $r(c)=\{v_0\}\in \CE $. 
Then $\{v\}\in H$ for each $v\in E^0$.
Consider the  quotient labeled space $(E,\CL,\CE/H)$.
\vskip .5pc 
\hskip 4pc 
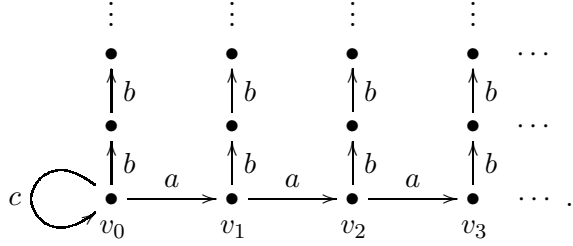
\xy /r0.38pc/:
(6,0)*+{\cdots\ . };(5,6)*+{\cdots};(5,12)*+{\cdots};
(-30,0)*+{\bullet}="V-3"; (-20,0)*+{\bullet}="V-2";
(-10,0)*+{\bullet}="V-1"; (0,0)*+{\bullet}="V0";
 (-20,6)*+{\bullet}="W-2"; (-30,6)*+{\bullet}="W-3";
(-10,6)*+{\bullet}="W-1"; (0,6)*+{\bullet}="W0";
(-20,12)*+{\bullet}="U-2";(-30,12)*+{\bullet}="U-3";
(-10,12)*+{\bullet}="U-1"; (0,12)*+{\bullet}="U0";
(0,16)*+{\vdots};(-10,16)*+{\vdots};(-20,16)*+{\vdots};(-30,16)*+{\vdots};
"V-3";"V-3"**\crv{(-30,0)&(-34,4)&(-38,0)&(-34,-4)&(-30,0)};
 ?>*\dir{>}\POS?(.5)*+!D{};
 "V-3";"V-2"**\crv{(-30,0)&(-20,0)};
 ?>*\dir{>}\POS?(.5)*+!D{};
 "V-2";"V-1"**\crv{(-20,0)&(-10,0)};
 ?>*\dir{>}\POS?(.5)*+!D{}; 
 "V-3";"W-3"**\crv{(-30,0)&(-30,6)};   
 ?>*\dir{>}\POS?(.5)*+!D{}; 
   "V-2";"W-2"**\crv{(-20,0)&(-20,6)};   
 ?>*\dir{>}\POS?(.5)*+!D{};  
 "V-1";"V0"**\crv{(-10,0)&(0,0)}; 
 ?>*\dir{>}\POS?(.5)*+!D{}; 
 "V-1";"W-1"**\crv{(-10,0)&(-10,6)}; 
 ?>*\dir{>}\POS?(.5)*+!D{};  
"V0";"W0"**\crv{(0,0)&(0,10)}; 
  ?>*\dir{>}\POS?(.5)*+!D{}; 
  "W-3";"U-3"**\crv{(-30,6)&(-30,12)};   
 ?>*\dir{>}\POS?(.5)*+!D{}; 
   "W-2";"U-2"**\crv{(-20,6)&(-20,12)};   
 ?>*\dir{>}\POS?(.5)*+!D{}; 
   "W-1";"U-1"**\crv{(-10,12)&(-10,12)};   
 ?>*\dir{>}\POS?(.5)*+!D{};  
   "W0";"U0"**\crv{(0,6)&(0,12)};   
 ?>*\dir{>}\POS?(.5)*+!D{};  
 (-38,0)*+{c};
 (-25,1.5)*+{a};(-15,1.5)*+{a};(-5,1.5)*+{a};
 (0.1,-2.5)*+{v_{3}}; (-9.9,-2.5)*+{v_{2}}; 
(-19.9,-2.5)*+{v_{1}};(-29.9,-2.5)*+{v_{0}}; 
(-28.5,3)*+{b};(-18.5,3)*+{b};(-8.5,3)*+{b};(1.5,3)*+{b};
(-28.5,9)*+{b};(-18.5,9)*+{b};(-8.5,9)*+{b};(1.5,9)*+{b};
\endxy
\vskip 1pc 
\noindent
Then  $[r(a)]=[r(a^2)]\neq [\emptyset]$, and 
since  
$$[r(a)] \cap r([r(a)],a)=[r(a)] \cap [r(a^2)] = [r(a)],$$ 
 we see that  $(a,[r(a)])$ is a loop. 
 On the other hand, for $b\in \CL_H([r(a)]E^{1})$, 
 $[r(b)]\neq[\emptyset]$ in  $\CE/H$, 
it follows from  $r([r(a)],b)=[r(ab)]\neq [\emptyset]$  
that the loop $(a,[r(a)])$ has an exit $b$ of type (i) 
 of Definition~\ref{qloop}. 
\vskip 1pc 
\noindent

\item[(ii)] 
 Let $H$ be the smallest hereditary saturated subset 
containing $r(c)=\{v_0\}$ in the following labeled space. 

\vskip 1pc 
\hskip 4pc 
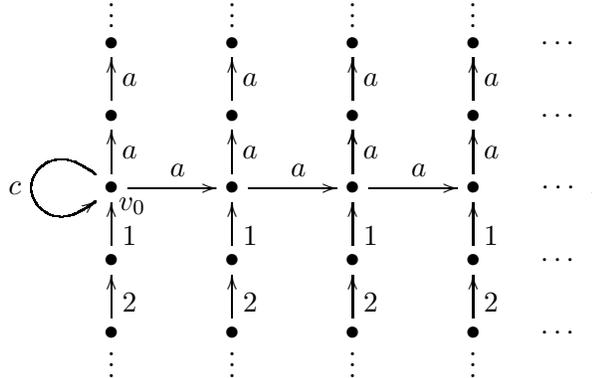
\xy /r0.38pc/:
(8,0)*+{\cdots \ . };
(7,6)*+{\cdots};
(7,12)*+{\cdots};
(7,-6)*+{\cdots};
(7,-12)*+{\cdots};
(-30,0)*+{\bullet}="V-3"; (-20,0)*+{\bullet}="V-2";
(-10,0)*+{\bullet}="V-1"; (0,0)*+{\bullet}="V0";
(-20,6)*+{\bullet}="W-2"; (-30,6)*+{\bullet}="W-3";
(-10,6)*+{\bullet}="W-1"; (0,6)*+{\bullet}="W0"; 
(-20,12)*+{\bullet}="U-2";(-30,12)*+{\bullet}="U-3";
(-10,12)*+{\bullet}="U-1"; (0,12)*+{\bullet}="U0";
(-20,-6)*+{\bullet}="Z-2";(-30,-6)*+{\bullet}="Z-3";
(-10,-6)*+{\bullet}="Z-1"; (0,-6)*+{\bullet}="Z0";
(-20,-12)*+{\bullet}="X-2";(-30,-12)*+{\bullet}="X-3";
(-10,-12)*+{\bullet}="X-1"; (0,-12)*+{\bullet}="X0";
(0,15)*+{\vdots};(-10,15)*+{\vdots};(-20,15)*+{\vdots};(-30,15)*+{\vdots};
(-30,-14)*+{\vdots};(-20,-14)*+{\vdots};(-10,-14)*+{\vdots};(0,-14)*+{\vdots};
"V-3";"V-3"**\crv{(-30,0)&(-34,4)&(-38,0)&(-34,-4)&(-30,0)};
 ?>*\dir{>}\POS?(.5)*+!D{};
 "V-3";"V-2"**\crv{(-30,0)&(-20,0)};
 ?>*\dir{>}\POS?(.5)*+!D{};
 "V-2";"V-1"**\crv{(-20,0)&(-10,0)};
 ?>*\dir{>}\POS?(.5)*+!D{}; 
 "V-3";"W-3"**\crv{(-30,0)&(-30,10)};   
 ?>*\dir{>}\POS?(.5)*+!D{}; 
   "V-2";"W-2"**\crv{(-20,0)&(-20,10)};   
 ?>*\dir{>}\POS?(.5)*+!D{};  
 "V-1";"V0"**\crv{(-10,0)&(0,0)}; 
 ?>*\dir{>}\POS?(.5)*+!D{}; 
 "V-1";"W-1"**\crv{(-10,0)&(-10,10)}; 
 ?>*\dir{>}\POS?(.5)*+!D{};  
  ?>*\dir{>}\POS?(.5)*+!D{}; 
 "V0";"W0"**\crv{(0,0)&(0,10)};  ?>*\dir{>}\POS?(.5)*+!D{}; 
 "W-3";"U-3"**\crv{(-30,8)&(-30,16)}; ?>*\dir{>}\POS?(.5)*+!D{}; 
 "W-2";"U-2"**\crv{(-20,8)&(-20,16)}; ?>*\dir{>}\POS?(.5)*+!D{}; 
 "W-1";"U-1"**\crv{(-10,8)&(-10,16)}; ?>*\dir{>}\POS?(.5)*+!D{};  
 "W0";"U0"**\crv{(0,8)&(0,16)}; ?>*\dir{>}\POS?(.5)*+!D{};  
 "Z-3";"V-3"**\crv{(-30,-8)&(-30,0)}; ?>*\dir{>}\POS?(.5)*+!D{};
 "Z-2";"V-2"**\crv{(-20,-8)&(-20,0)}; ?>*\dir{>}\POS?(.5)*+!D{}; 
 "Z-1";"V-1"**\crv{(-10,-8)&(-10,0)}; ?>*\dir{>}\POS?(.5)*+!D{};  
 "Z0";"V0"**\crv{(0,-8)&(0,0)}; ?>*\dir{>}\POS?(.5)*+!D{};  
 "X-3";"Z-3"**\crv{(-30,-16)&(-30,-8)}; ?>*\dir{>}\POS?(.5)*+!D{};
 "X-2";"Z-2"**\crv{(-20,-16)&(-20,-8)}; ?>*\dir{>}\POS?(.5)*+!D{}; 
 "X-1";"Z-1"**\crv{(-10,-16)&(-10,-8)}; ?>*\dir{>}\POS?(.5)*+!D{};  
 "X0";"Z0"**\crv{(0,-16)&(0,-8)}; ?>*\dir{>}\POS?(.5)*+!D{};  
(-24.5,1.5)*+{a};(-14.5,1.5)*+{a};(-4.5,1.5)*+{a};
(-38,0)*+{c};
(-28.3,-1.5)*+{v_{0}};
(-28.5,3)*+{a};(-18.5,3)*+{a};(-8.5,3)*+{a};(1.5,3)*+{a};
(-28.5,9)*+{a};(-18.5,9)*+{a};(-8.5,9)*+{a};(1.5,9)*+{a};
(-28.5,-4)*+{1};(-18.5,-4)*+{1};(-8.5,-4)*+{1};(1.5,-4)*+{1};
(-28.5,-9.5)*+{2};(-18.5,-9.5)*+{2};(-8.5,-9.5)*+{2};(1.5,-9.5)*+{2};
\endxy
\vskip 1pc 
\noindent
Then $(a, [r(1)])$ is a loop 
in the quotient labeled space $(E,\CL,\CE/H)$ 
such that $[r(1)] \subsetneq  r([r(1)],a)=[(r(1a)]$. 
That is, the loop $(a, [r(1)])$  has 
an exit of type (ii) of  Definition \ref{qloop}.
\end{enumerate}
\end{ex}

\vskip 1pc 
 
If a directed  graph $E$ has a loop with an exit, 
the exit gives rise to an infinite projection in 
the graph $C^*$-algebra $C^*(E)$. 
The same is true for labeled graph $C^*$-algebras 
(\cite[Proposition 3.5]{JKK}), and more generally for 
quotient lebeled graph $C^*$-algebras as we see from the 
following proposition.  

\vskip 1pc 

\begin{prop} \label{QLE} 
Let  $(E,\CL,\CE/R)$ be a  quotient labeled space.  
If $(\af,[A])$ is a loop with an exit, then  
$p_{[A]}$ is an infinite projection in $C^*(E,\CL,\CE/R)$.
In particular, if  $(E,\CL,\CE/R)$ is disagreeable, 
every loop  $(\af,[A])$ has an exit and 
$p_{[A]}$ is an infinite projection.
\end{prop}

\begin{proof} 
If $(\af,[A])$ has an exit of type (i) with 
a path $\bt \in \CL_R([A]E^{|\af|})$ such that 
$\bt \neq \af$  and $r([A],\bt) \neq [\emptyset]$,  
then  the projection $p_{[A]}$ is infinite since
$$p_{r([A],\af)} \geq\, p_{[A]} 
=  \sum_{|\dt|=|\af|}s_{\dt}p_{r([A],\dt)}s_{\dt}^* 
> s_{\af}p_{r([A],\af)}s_{\af}^* \sim p_{r([A],\af)}.$$
If $(\af,[A])$ has an exit of type (ii), that is 
$[A] \subsetneq r([A], \af)$, then 
$$p_{r([A], \af)}\sim s_{\af}p_{r([A], \af)}s_{\af}^* 
\leq p_{[A]}\lneq p_{r([A], \af)},$$ 
which shows that $p_{[A]}$ is infinite. 

To prove the second assertion, 
let $(E,\CL,\CE)$ be a disagreeable labeled space. 
If $A\subsetneq r([A],\af)$, $\af$ has an exit of type (ii). 
Assume that $[A]= r([A],\af)$.
Then by Proposition~\ref{prop-disagreeable} 
(see Remark~\ref{q-disagreeable}), there is an $n\geq 1$ 
such that 
$\CL([A]E^{|\af|n})\neq \{\af^n\}$. 
Let $m$ be the smallest integer among such $n$'s with 
$\CL([A]E^{|\af|n})\neq \{\af^n\}$, 
and choose  a path $\bt\in \CL([A]E^{|\af|m})$ with $\bt\neq \af^m$. 
If $m=1$, or $m>1$ and $\bt_{[1,|\af|]}\neq \af$, 
then $\bt_{[1,|\af|]}$ is an exit of $(\af,[A])$. 
If $m>1$ and $\bt=\af\bt''$ for some $\bt''$, 
then from $[A]= r([A],\af)$ 
we see that $\bt'' \in \CL(AE^{|\af|(m-1)})$.  
Then $\bt''\neq \af^{m-1}$, that is 
$\CL([A]E^{|\af|(m-1)})\neq \{\af^{m-1}\}$ follows, 
but this contradicts to the choice of $m$. 
\end{proof}

\vskip 1pc
 
\subsection{Labeled graph $C^*$-algebras with 
the property (IH)} 

As mentioned in Introduction, it is well known that 
for a graph $E$ satisfying condition (L), 
the $C^*$-algebra 
$C^*(E)$ has the property (IH) if and only if 
every vertex connects to a loop. 
We will prove in Theorem~\ref{thm-IH} 
that for a disagreeable labeled space $(E,\CL,\CE)$, 
if every vertex connects to a loop in the sense of 
Definition~\ref{VCL} below, then 
$C^*(E,\CL,\CE)$ has the property (IH).  
Its converse will be shown to be true 
under an additional condition. 
 To prove Theorem~\ref{thm-IH}, 
we need the following useful proposition 
which can be obtained by similar arguments 
used in \cite[Theorem 6.9]{BP2} and 
\cite[Proposition 5.3]{BPRS} with slight modifications. 
We provide a detailed proof here for the reader's convenience.

\vskip 1pc

\begin{prop}\label{DHP} Let $(E,\CL,\CE )$ be a  disagreeable labeled space.
Then every nonzero hereditary $C^*$-subalgebra of $C^*(E,\CL,\CE )$ 
contains a nonzero projection 
$p$ such that $0\neq s_\mu p_A s_\mu^* \preceq p$ for some 
$\mu\in \CL^*(E)$ and $A\in \CE$ with $A\subset r(\mu)$.
\end{prop}
 
\begin{proof}  Let $B$ be a nonzero hereditary $C^*$-subalgebra  of  
$C^*(E,\CL,\CE )$ and fix a positive element $a\in B$ with $\|\Phi(a)\| = 1$. 
Choose a positive element 
$b \in {\rm span}\{s_{\af}p_{A}s_{\bt}^* :~ \af, \bt \in \CL^*(E)~
\text{and}~ A \subseteq r(\af)\cap r(\bt) \}$ 
so that $\|a-b\|<\frac{1}{4}$. 
From \cite[Proposition 2.4.(ii) and (iii)]{BP2}, 
we may write 
$b=\sum_{(\af,[w]_l,\bt) \in F}c_{(\af,[w]_l,\bt)}s_{\af}p_{[w]_l}s_{\bt}^*$, 
where $F$ is a finite subset of 
$\CL^*(E) \times \Omega_l \times \CL^*(E)$ for some $l \geq 1$. 
Let $b_0= \Phi(b) \geq 0 $. 
Since $\Phi$ is norm-decreasing, 
we have 
$$|1- \| b_0 \|| = |\| \Phi(a) \|-\| \Phi(b)\| | \leq \| \Phi(a-b) \| \leq \| a-b\|
 < \frac{1}{4},$$
and hence $\| b_0 \| \geq \frac{3}{4} $.
 Let  $k = \max \{ |\af|, |\bt| :(\af,[w]_l,\bt) \in F \}$.
 Applying Definition \ref{def-representation}.(iv) and 
 changing $F$ (if necessary),
 we can choose a $k \in \mathbb{N}$ so that 
  $\min \{ |\af|, |\bt|\}=k$   for every $(\af,[w]_l,\bt) \in F$. 
 Let  $M=\max \{ |\af|, |\bt| :(\af,[w]_l,\bt) \in F \}.$ 
 Applying \cite[Proposition 2.4.(iii)]{BP2} again, 
we may choose  $m \geq M$  large enough so that
$$b_0 \in \oplus_{ \{w: (\af,[w]_l,\bt) \in F\}} \mathcal{F}^k([w]_m) $$ 
(see (\ref{Fk}) for $\mathcal{F}^k([w]_m)$).
Now, $\| b_0 \|$ must be attained in some summand $\mathcal{F}^k([v]_m)$. 
Let $b_1$ be the component of $b_0$ in  $\mathcal{F}^k([v]_m)$ so that  
$\| b_0 \|=\| b_1 \|$ and note that $b_1 \geq 0$. 
Then we can choose a projection $r \in C^*(b_1) \subseteq \mathcal{F}^k([v]_m)$ 
such that $rb_1r=\| b_1 \| r$.
 Since $b_1$ is a finite sum of $s_{\af}p_{[v]_m}s_{\bt}^*$, 
 we can write $r$ as a sum $\sum c_{\af\bt}s_{\af}p_{[v]_m}s_{\bt}^*$ 
 over all pairs of paths in 
$$G=\{\af \in \CL(E^k) : ~\text{either} ~ (\af, [v]_m, \bt) \in F 
~\text{or}~ (\bt, [v]_m, \af) \in F \}.$$  
Note that $rb_0r=rb_1r$ and the $G \times G$-matrix $(c_{\af\bt})$ is also 
a projection in a finite dimensional matrix algebra 
$\CF^k([v]_m)={\rm span} \{s_{\af}p_{[v]_m}s_{\bt}^* : \af, \bt \in G\}$.

Since $[v]_m $ is disagreeable, we may choose a path  $\ld \in \CL^*(E)$   
with $|\ld|>M$ so that $\ld$  has no factorization $\ld=\ld'\ld''=\ld''\dt$ 
for some $|\ld'|,|\dt| \leq m$.
  Then because 
  ${\rm span}\{s_{\af\ld}p_{r([v]_m,\ld)} s_{\bt\ld}^*: \af,\bt \in G \}$ 
  is also   a finite dimensional matrix algebra generated by  
  the family of non-zero matrix units 
  $\{s_{\af\ld}p_{r([v]_m,\ld)} s_{\bt\ld}^*: \af,\bt \in G\}$, 
  $$Q=\sum_{\af,\bt \in G}c_{\af\bt}s_{\af\ld}p_{r([v]_m,\ld)}s_{\bt\ld}^*$$ 
is a projection satisfying 
$$ r = \sum  c_{\af\bt}s_{\af}p_{[v]_m}s_{\bt}^* 
 =\sum c_{\af\bt}s_{\af}(s_{\ld}p_{r([v]_m, \ld)}s_{\ld}^* +(p_{[v]_m} 
 -s_{\ld}p_{r([v]_m, \ld)}s_{\ld}^*))s_{\bt}^* \geq Q.$$

  We claim that for $(\mu,[v]_m,\nu)\in F$,
  $$Qs_{\mu}p_{[v]_m}s_{\nu}^*Q=0 ~\text{unless}~ 
  |\mu|=|\nu|=k ~\text{and}~ [v]_m \subseteq r(\mu) \cap r(\nu).$$ 
Suppose that $(\mu, [v]_m, \nu) \in F$ with $|\mu|\neq|\nu|$. 
We may assume $|\mu|=k$ because either $\mu$ or $\nu$ has length k. 
Since $s_{\bt\ld}^*s_{\mu} \neq 0$ if and only if $\bt=\mu$, we have
\begin{align*}
Q(s_{\mu}p_{[v]_m}s_{\nu}^*)Q 
  &=
  (\sum_{\af',\bt' \in G}c_{\af'\bt'}s_{\af'\ld}p_{r([v]_m,\ld)}s_{\bt'\ld}^*)
  (s_{\mu}p_{[v]_m}s_{\nu}^*)(\sum_{\af,\bt \in G}c_{\af\bt}s_{\af\ld}p_{r([v]_m,\ld)}s_{\bt\ld}^*)\\ 
  &=  (\sum_{\af' \in G}c_{\af'\mu}s_{\af'\ld}p_{r([v]_m,\ld)}s_{\mu\ld}^*s_{\mu}p_{[v]_m}s_{\nu}^*)
  (\sum_{\af,\bt \in G}c_{\af\bt}s_{\af\ld}p_{r([v]_m,\ld)}s_{\bt\ld}^*)\\ 
  &=\sum_{\af,\bt \in G}c_{\af\bt}(\sum_{\af' \in G}c_{\af'\mu}s_{\af'\ld}p_{r([v]_m,\ld)}
  s_{\nu\ld}^*)s_{\af\ld}p_{r([v]_m,\ld)}s_{\bt\ld}^*.  
\end{align*}
To be $s_{\nu\ld}^*s_{\af\ld} \neq 0$, $\nu\ld$ must extend $\af\ld$, so that
$\nu\ld=\af\ld\dt$ for some $\dt \in \CL^*(E)$. 
On the other hand, we may say that $\nu=\af\ld'$ where $\ld=\ld'\ld''$ 
for some $\ld',\ld'' \in \CL^*(E)$ since $|\nu|>|\af|=k$. 
As $\af\ld'\ld=\nu\ld=\af\ld\dt=\af\ld'\ld''\dt$, 
we have 
$$\ld=\ld'\ld''=\ld''\dt$$ with $|\ld'|=|\dt|$. 
Because $|\nu|=|\af\ld'|\leq M $ with $|\af|=k$, 
we know $|\ld'| \leq M-k \leq m$,
which contradicts to the fact that $\ld$ is disagreeable  for $[v]_m$.
      
Thus, we see that
 $$ QbQ=Qb_1Q=Qrb_1rQ=\| b_1\| rQ=\| b_0 \| Q \geq \frac{3}{4}Q.$$
Since $\|a-b\| <\frac{1}{4} $, 
we have $QaQ \geq QbQ -\frac{1}{4}Q \geq \frac{1}{2}Q$. 
This implies that $QaQ$  is invertible in $QC^*(E,\CL,\CE )Q$.
Let $c$ be the inverse of $QaQ$ in $QC^*(E,\CL,\CE )Q$ and 
put $v=c^{\frac{1}{2}}Qa^{\frac{1}{2}}$. 
Then $v^*v=a^{\frac{1}{2}}QcQa^{\frac{1}{2}} \leq \|c \| a$, 
and hence $v^*v \in B$.
Since $$v^*v \sim vv^*=c^{\frac{1}{2}}QaQc^{\frac{1}{2}}=Q,$$ 
the hereditary $C^*$-subalgebra $B$ contains a non-zero projection equivalent to $Q$.
Note that $Q$ belongs to 
the finite dimensional subalgebra 
$C:={\rm span}\{s_{\af\ld}p_{r([v]_m,\ld)} s_{\bt\ld}^*: \af,\bt \in G \}$
for which the elements $\{s_{\af\ld}p_{r([v]_m,\ld)} s_{\bt\ld}^*\}$ forms 
a matrix unit. 
This means that $Q$ dominates a minimal projection in $C$. 
Since every minimal projection in $C$ is  equivalent to a minimal projection of the form 
$s_{\af\ld}p_{r([v]_m,\ld)} s_{\af\ld}^*$,  
the hereditary subalgebra $B$ also contains a projection 
equivalent to the desired form.
\end{proof}

\vskip 1pc

\begin{prop}\label{simple PI} 
Let $(E,\CL,\CE)$ be a disagreeable and strongly cofinal labeled space.
Then the following are equivalent: 
\begin{enumerate}
\item[(a)] $C^*(E,\CL,\CE )$ is purely infinite.
\item[(b)] Every nonzero projection $p_A$ is infinite 
for $A\in \CE$.
\end{enumerate} 
\end{prop}
 
\begin{proof} 
We only need to show that $(b)$ implies (a). 
By Proposition~\ref{DHP}, every nonzero hereditary subalgebra contains a projection 
$p$ such that $p_A\preceq p$ for a nonempty set $A\in \CE$. 
Thus $p$ should be infinite, and 
$C^*(E,\CL,\CE )$ has the property (IH). 
Since  $C^*(E,\CL,\CE )$ is simple, it is purely infinite.
\end{proof}

\vskip 1pc

\begin{dfn}\label{VCL} 
Let  $(E,\CL, \CE)$ be a labeled space.
We say that {\it every vertex connects to a loop} 
in $(E,\CL, \CE)$ if for each $[v]_l$, 
there exist  a loop $(\af,A)$ and 
a finitely many  paths 
$\dt_1, \dots, \dt_m \in \CL^*(E)$ 
such that $\dt_i$ is not an initial path of $\dt_j$ for 
$j\neq i$ and 
\begin{eqnarray}\label{vcl}
A \subseteq \cup_{i=1}^m r([v]_l,\dt_i).
\end{eqnarray}

\noindent 
This definition naturally extends to 
quotient labeled spaces $(E,\CL,\CE/R)$.
\end{dfn}

\vskip 1pc
\noindent 
Note that the condition in the above definition 
for a path  $\dt_i$  to be 
an initial path of none of the rest of the paths implies that 
 $s_{\dt_i}^* s_{\dt_j}=0$ for $i\neq j$ in $C^*(E,\CL,\CE)$. 
This is the case if $\dt_i$'s are distinct paths with the same length.

\vskip 1pc

\begin{remark}\label{rmk-vcl} 
 In \cite[Definition 6.6]{BP2},    
it is said that every vertex connects to a repeatable path if   
for each $[v]_m$, there exist $w \in E^0$, $L \geq 1$, and 
 $\af, \dt \in \CL^*(E)$ such that 
\begin{eqnarray}\label{connect to rp} 
[w]_L \subseteq  r([v]_m, \dt)
\ \text{and }\  (\af,[w]_l)\ \text{is a loop  
 for}\ l \geq L.
\end{eqnarray}
This property is obviously stronger than the one 
(every vertex connects to a loop)
given in Definition~\ref{VCL}, and 
actually strictly stronger as we will see in Example~\ref{ex-VCL}.  
It is known in \cite[Theorem 6.9]{BP2} that 
if a labeled space $(E,\CL,\CE )$ is disagreeable, 
strongly cofinal, 
and every vertex connects to a repeatable path, then 
the simple $C^*$-algebra $C^*(E,\CL,\CE )$ has the property  (IH).
We will extend this result to 
disagreeable labeled spaces in which 
every vertex connects to a loop (see Theorem~\ref{thm-IH}). 
\end{remark}

\vskip 1pc

\begin{ex}\label{ex-VCL} 
In the following labeled space $(E,\CL, \CE)$,  
 if $(\af,[w]_l)$ is a loop, then 
$w=v_k$ and $\af\in \{a,b\}^n$ for some $k,n \geq 1$. 

\vskip 1pc \hskip 1pc
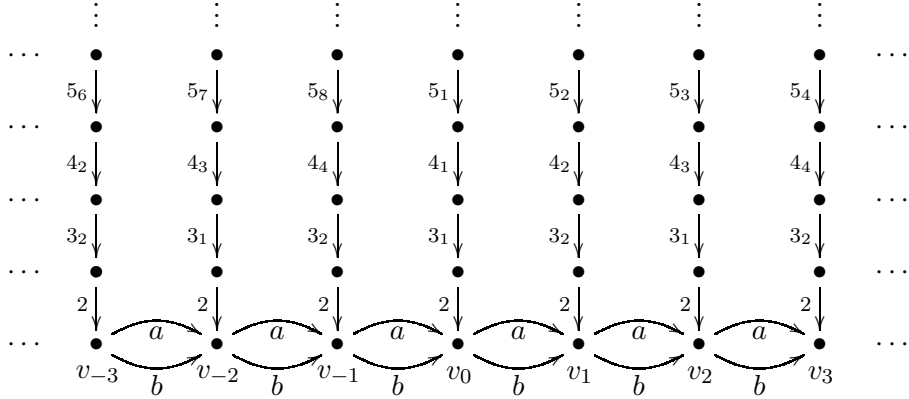
\xy /r0.38pc/:
(-36,0)*+{\cdots};(36,0)*+{\cdots};
(36,6)*+{\cdots};(-36,6)*+{\cdots};
(-36,12)*+{\cdots};(36,12)*+{\cdots};
(-36,18)*+{\cdots};(36,18)*+{\cdots};
(36,24)*+{\cdots};(-36,24)*+{\cdots};
(0,28)*+{\vdots};(10,28)*+{\vdots};
(20,28)*+{\vdots};(30,28)*+{\vdots};
(-10,28)*+{\vdots};(-20,28)*+{\vdots}; 
(-30,28)*+{\vdots};
(-30,0)*+{\bullet}="V-3"; 
(-20,0)*+{\bullet}="V-2";
(-10,0)*+{\bullet}="V-1"; 
(0,0)*+{\bullet}="V0";
(10,0)*+{\bullet}="V1"; 
(20,0)*+{\bullet}="V2";
(30,0)*+{\bullet}="V3";
(-30,6)*+{\bullet}="V-31"; 
(-20,6)*+{\bullet}="V-21";
(-10,6)*+{\bullet}="V-11"; 
(0,6)*+{\bullet}="V01";
(10,6)*+{\bullet}="V11"; 
(20,6)*+{\bullet}="V21";
(30,6)*+{\bullet}="V31";
(-30,12)*+{\bullet}="V-32"; 
(-20,12)*+{\bullet}="V-22";
(-10,12)*+{\bullet}="V-12"; 
(0,12)*+{\bullet}="V02";
(10,12)*+{\bullet}="V12"; 
(20,12)*+{\bullet}="V22";
(30,12)*+{\bullet}="V32";
(0,18)*+{\bullet}="V03";
(-10,18)*+{\bullet}="V-13";
(10,18)*+{\bullet}="V13";
(-20,18)*+{\bullet}="V-23";
(20,18)*+{\bullet}="V23";
(-30,18)*+{\bullet}="V-33";
(30,18)*+{\bullet}="V33";
(0,24)*+{\bullet}="V04";
(-30,6)*+{\bullet}="V-31";
(-10,24)*+{\bullet}="V-14";
(10,24)*+{\bullet}="V14";
(-20,24)*+{\bullet}="V-24";
(20,24)*+{\bullet}="V24";
(-30,24)*+{\bullet}="V-34";
(30,24)*+{\bullet}="V34";
"V03";"V02"**\crv{(0,18)&(0,6)}; ?>*\dir{>}\POS?(.5)*+!D{}; 
 "V13";"V12"**\crv{(10,18)&(10,6)}; ?>*\dir{>}\POS?(.5)*+!D{}; 
 "V23";"V22"**\crv{(20,18)&(20,6)}; ?>*\dir{>}\POS?(.5)*+!D{}; 
 "V33";"V32"**\crv{(30,18)&(30,6)}; ?>*\dir{>}\POS?(.5)*+!D{}; 
 "V-13";"V-12"**\crv{(-10,18)&(-10,6)}; ?>*\dir{>}\POS?(.5)*+!D{}; 
 "V-23";"V-22"**\crv{(-20,18)&(-20,6)}; ?>*\dir{>}\POS?(.5)*+!D{}; 
 "V-33";"V-32"**\crv{(-30,18)&(-30,6)}; ?>*\dir{>}\POS?(.5)*+!D{};  
 "V04";"V03"**\crv{(0,24)&(0,24)}; ?>*\dir{>}\POS?(.5)*+!D{}; 
 "V14";"V13"**\crv{(10,24)&(10,24)}; ?>*\dir{>}\POS?(.5)*+!D{}; 
 "V24";"V23"**\crv{(20,24)&(20,24)}; ?>*\dir{>}\POS?(.5)*+!D{}; 
 "V34";"V33"**\crv{(30,24)&(30,24)}; ?>*\dir{>}\POS?(.5)*+!D{};  
 "V-14";"V-13"**\crv{(-10,24)&(-10,12)}; ?>*\dir{>}\POS?(.5)*+!D{};
 "V-24";"V-23"**\crv{(-20,24)&(-20,12)}; ?>*\dir{>}\POS?(.5)*+!D{}; 
 "V-34";"V-33"**\crv{(-30,24)&(-30,12)}; ?>*\dir{>}\POS?(.5)*+!D{};       
 "V-32";"V-31"**\crv{(-30,12)&(-30,6)}; ?>*\dir{>}\POS?(.5)*+!D{}; 
 "V-31";"V-3"**\crv{(-30,6)&(-30,0)}; ?>*\dir{>}\POS?(.5)*+!D{}; 
 "V-22";"V-21"**\crv{(-20,12)&(-20,6)}; ?>*\dir{>}\POS?(.5)*+!D{};
 "V-21";"V-2"**\crv{(-20,6)&(-20,0)}; ?>*\dir{>}\POS?(.5)*+!D{};
 "V-12";"V-11"**\crv{(-10,12)&(-10,6)}; ?>*\dir{>}\POS?(.5)*+!D{}; 
 "V-11";"V-1"**\crv{(-10,6)&(-10,0)}; ?>*\dir{>}\POS?(.5)*+!D{};
 "V02";"V01"**\crv{(0,12)&(0,6)}; ?>*\dir{>}\POS?(.5)*+!D{}; 
 "V01";"V0"**\crv{(0,6)&(0,0)}; ?>*\dir{>}\POS?(.5)*+!D{};
 "V12";"V11"**\crv{(10,12)&(10,6)}; ?>*\dir{>}\POS?(.5)*+!D{}; 
 "V11";"V1"**\crv{(10,6)&(10,0)}; ?>*\dir{>}\POS?(.5)*+!D{};
 "V22";"V21"**\crv{(20,12)&(20,6)}; ?>*\dir{>}\POS?(.5)*+!D{}; 
 "V21";"V2"**\crv{(20,6)&(20,0)}; ?>*\dir{>}\POS?(.5)*+!D{};
 "V32";"V31"**\crv{(30,12)&(30,6)}; ?>*\dir{>}\POS?(.5)*+!D{}; 
 "V31";"V3"**\crv{(30,6)&(30,0)}; ?>*\dir{>}\POS?(.5)*+!D{};  
 "V-3";"V-2"**\crv{(-30,0)&(-25,4)&(-20,0)}; ?>*\dir{>}\POS?(.5)*+!D{};
 "V-2";"V-1"**\crv{(-20,0)&(-15,4)&(-10,0)}; ?>*\dir{>}\POS?(.5)*+!D{};
 "V-1";"V0"**\crv{(-10,0)&(-5,4)&(0,0)}; ?>*\dir{>}\POS?(.5)*+!D{}; 
 "V0";"V1"**\crv{(0,0)&(5,4)&(10,0)}; ?>*\dir{>}\POS?(.5)*+!D{};
 "V1";"V2"**\crv{(10,0)&(15,4)&(20,0)}; ?>*\dir{>}\POS?(.5)*+!D{};
 "V2";"V3"**\crv{(20,0)&(25,4)&(30,0)}; ?>*\dir{>}\POS?(.5)*+!D{};  
 "V-3";"V-2"**\crv{(-30,0)&(-25,-4)&(-20,0)}; ?>*\dir{>}\POS?(.5)*+!D{};
 "V-2";"V-1"**\crv{(-20,0)&(-15,-4)&(-10,0)}; ?>*\dir{>}\POS?(.5)*+!D{};
 "V-1";"V0"**\crv{(-10,0)&(-5,-4)&(0,0)}; ?>*\dir{>}\POS?(.5)*+!D{}; 
 "V0";"V1"**\crv{(0,0)&(5,-4)&(10,0)}; ?>*\dir{>}\POS?(.5)*+!D{};
 "V1";"V2"**\crv{(10,0)&(15,-4)&(20,0)}; ?>*\dir{>}\POS?(.5)*+!D{};
 "V2";"V3"**\crv{(20,0)&(25,-4)&(30,0)}; ?>*\dir{>}\POS?(.5)*+!D{};   
 ;(-25,1)*+{a};(-15,1)*+{a};(-5,1)*+{a};(5,1)*+{a};
 (15,1)*+{a};(25,1)*+{a};
(-25,-3.3)*+{b};(-15,-3.3)*+{b};(-5,-3.3)*+{b};
(5,-3.3)*+{b};(15,-3.3)*+{b};(25,-3.3)*+{b};
(0.1,-2.5)*+{v_0}; (10.1,-2.5)*+{v_1};(-9.9,-2.5)*+{v_{-1}};
(-19.9,-2.5)*+{v_{-2}};(-29.9,-2.5)*+{v_{-3}};
(20.1,-2.5)*+{v_{2}};(30.1,-2.5)*+{v_{3}};    
(-31.1,3.4)*+{{}_{2}};(-21.1,3.4)*+{{}_{2}};(-11.1,3.4)*+{{}_{2}};
(-1.1,3.4)*+{{}_{2}};(8.9,3.4)*+{{}_{2}}; (18.9,3.4)*+{{}_{2}};
(28.9,3.4)*+{{}_{2}};
(-31.5,9)*+{{}_{3_2}};(-21.5,9)*+{{}_{3_1}};
(-11.5,9)*+{{}_{3_2}};(-1.5,9)*+{{}_{3_1}}; (8.5,9)*+{{}_{3_2}};
(18.5,9)*+{{}_{3_1}}; (28.5,9)*+{{}_{3_2}}; 
(-31.5,15)*+{{}_{4_2}};(-21.5,15)*+{{}_{4_3}}; (-11.5,15)*+{{}_{4_4}};
(-1.5,15)*+{{}_{4_1}}; (8.5,15)*+{{}_{4_2}};
(18.5,15)*+{{}_{4_3}};(28.5,15)*+{{}_{4_4}};
(-31.5,21)*+{{}_{5_6}};(-21.5,21)*+{{}_{5_7}};
(-11.5,21)*+{{}_{5_8}}; (-1.5,21)*+{{}_{5_1}};(8.5,21)*+{{}_{5_2}};
(18.5,21)*+{{}_{5_3}};(28.5,21)*+{{}_{5_4}};
\endxy 
\vskip 1pc 
\noindent 
But for a fixed $n$, if $l'$ is large enough, 
then $(\af, [v_k]_{l'})$ fails to be a loop, 
accordingly there is no  vertex that connects to 
a repeatable path in the sense of (\ref{connect to rp}). 
On the other hand, 
every vertex connects to a loop 
$([v_j]_l, a^n)$ for  some $j,l,n\geq 1$.
Since  $(E,\CL, \CE)$ is disagreeable and strongly cofinal, 
the  $C^*$-algebra $C^*(E,\CL,\CE )$ is simple, 
and we will see that  it is purely infinite
by Theorem~\ref{thm-IH} below.
\end{ex}

\vskip 1pc
Recall that a graph $C^*$-algebra $C^*(E)$ 
has the property  (IH) if and only if 
$E$ satisfies condition (L) and 
every vertex connects to a loop. 
 In the following theorem we show that 
a labeled graph $C^*$-algebra $C^*(E,\CL,\CE)$ has (IH) if 
$(E,\CL,\CE)$ is disagreeable and 
every vertex connects to a loop.
For the converse, we need the following additional assumptions 
for $(E,\CL,\CE)$:
$$ 
(*)\left\{\begin{array}{l}
\text{if $A\in\CE$, $A\neq \emptyset$, 
   then  $B\subset A$ for some $B\in\CE_{\rm min}$, and}\\ 
\text{if $A\in \CE_{\rm min}$, then   
$r(A,\af)\in \big\{\cup_{i=1}^n B_i: B_i\in \CE_{\rm min},\, n\geq 1\big\}$ 
for all  $\af\in \CL^*(E)$.}
\end{array}\right.
$$
Note that the condition $(*)$  also 
makes sense for quotient labeled spaces $(E,\CL,\CE/R)$ 
(with $\CE/R$ and $\CL_R^*(E)$ in place of $\CE$ and $\CL^*(E)$, 
respectively). 
While $(*)$ seems restrictive,  it is 
not too restrictive to imply that every nonempty   
 $A\in \CE$ is a finite union of minimal sets. 
For example, in the following disagreeable labeled space,  

\vskip 1pc
\hskip 5pc 
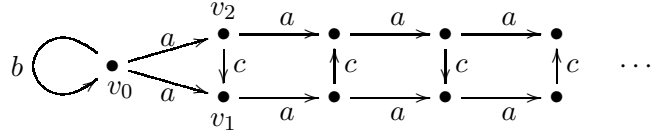
\xy /r0.5pc/:
 (20,0)*+{\bullet}="V1";
 (27,-2)*+{\bullet}="V2";(27,2)*+{\bullet}="V22"; 
 (34,-2)*+{\bullet}="V3";(34,2)*+{\bullet}="V32";
 (41,-2)*+{\bullet}="V4";(41,2)*+{\bullet}="V42";
 (48,-2)*+{\bullet}="V5";(48,2)*+{\bullet}="V52";
  "V1";"V1"**\crv{(20,0)&(17,3)&(14,0)&(17,-3)&(20,0)};?>*\dir{>}\POS?(.5)*+!D{};
 "V1";"V2"**\crv{(20,0)&(27,-2)};?>*\dir{>}\POS?(.5)*+!D{};
 "V2";"V3"**\crv{(27,-2)&(34,-2)};?>*\dir{>}\POS?(.5)*+!D{};
 "V3";"V4"**\crv{(34,-2)&(41,-2)};?>*\dir{>}\POS?(.5)*+!D{};
 "V4";"V5"**\crv{(41,-2)&(48,-2)};?>*\dir{>}\POS?(.5)*+!D{};
 "V1";"V22"**\crv{(20,0)&(27,2)};?>*\dir{>}\POS?(.5)*+!D{};
 "V22";"V2"**\crv{(27,2)&(27,0)};?>*\dir{>}\POS?(.5)*+!D{};
 "V22";"V32"**\crv{(27,2)&(34,2)};?>*\dir{>}\POS?(.5)*+!D{}; 
 "V3";"V32"**\crv{(34,-2)&(34,2)};?>*\dir{>}\POS?(.5)*+!D{};
 "V32";"V42"**\crv{(34,2)&(41,2)};?>*\dir{>}\POS?(.5)*+!D{};
 "V42";"V52"**\crv{(41,2)&(48,2)};?>*\dir{>}\POS?(.5)*+!D{};
 "V42";"V4"**\crv{(41,2)&(41,-2)};?>*\dir{>}\POS?(.5)*+!D{};
 "V5";"V52"**\crv{(48,-2)&(48,2)};?>*\dir{>}\POS?(.5)*+!D{};
  (14,0)*+{b},
  (23.5,-1.8)*+{a},(23.5,1.7)*+{a},(31,-3)*+{a},(31,3)*+{a},
  (38,-3)*+{a},(45,-3)*+{a},(38,3)*+{a},(45,3)*+{a},
  (28,0)*+{c},(35,0)*+{c},(42,0)*+{c},(49,0)*+{c},
  (20.5,-1.5)*+{v_0},(27,3.5)*+{v_2},(27,-3.5)*+{v_1}, 
  (53,0)*+{\cdots},  
\endxy
\vskip 1pc 
\noindent
each vertex set  $\{v\}\in \CE$, $v\in E^0$, is 
a minimal set in $\CE$, hence 
every nonempty set $A\in \CE$ contains a minimal set and 
obviously every relative range with respect to a minimal set is a finite 
union of minimal sets. 
But not every $A\in \CE$, like $A=r(a)$, 
is a finite union of minimal sets. 
\vskip 2pc 

\begin{thm}\label{thm-IH} 
For a disagreeable labeled space $(E,\CL,\CE )$, 
we have the following:
\begin{enumerate}
\item[(a)]
If every vertex connects to a loop,  then   
$C^*(E,\CL,\CE )$ has the property (IH).
 Moreover every nonzero hereditary $C^*$-subalgebra 
 of $C^*(E,\CL,\CE )$ contains an infinite projection 
 equivalent to a projection $p_A$ for some $A\in \CE $.
\item[(b)] 
Assume further that  $(E,\CL,\CE )$ satisfies $(*)$. 
If $C^*(E,\CL,\CE )$ has the property (IH), then 
every vertex connects to a loop.
\end{enumerate}
\end{thm}
 
\begin{proof} 
(a) 
Let $C^*(E,\CL,\CE)=C^*(s_a, p_A)$. 
We first show that the projections $p_{[v]_l}$ of 
generalized vertices $[v]_l$ are all infinite.
Given a generalized vertex  $[v]_l$, choose a path $\gm\in \CL^l(E)$ 
such that $[v]_l\subset r(\gm)$. 
Since $[v]_l$ connects to a loop by our assumption, 
there exist a loop 
$(\af, A)$ and finitely many paths $\dt_1, \dots, \dt_m$  such that 
$s_{\dt_i}^*s_{\dt_j}=0$ whenever $i\neq j$. 
Set 
$$A_i:=A\cap r([v]_l,\delta_i)\subset A\cap r(\gm\delta_i)\ \text{ and } \
B_i:=A_i\setminus \cup_{j=1}^{i-1} A_j \ (B_1:=A_1)$$ 
for $i=1, \dots, m$.  
Then the projections $p_{B_i}$'s are mutually orthogonal and 
\begin{align*} 
p_A  &= \sum_i \, p_{B_i} 
      \,\leq \,\sum_i \, p_{B_i} s_{\gm \dt_i}^* s_{\gm \dt_i}p_{B_i} \\ 
     &=\, \Big(\sum_i p_{B_i}s_{\gm \dt_i}^* \Big) \Big(\sum_i p_{B_i}s_{\gm \dt_i}^*\Big)^*             
      \sim \Big(\sum_i p_{B_i}s_{\gm \dt_i}^*\Big)^* \Big(\sum_i p_{B_i}s_{\gm \dt_i}^*\Big) \\
     &=\, \sum_i \, s_{\gm \dt_i} p_{B_i} s_{\gm \dt_i}^*  
     =\,  s_{\gm}\cdot \sum_i\,  s_{\dt_i} p_{B_i} s_{\dt_i}^* \cdot s_{\gm}^* \\     
     &  \leq  \, s_{\gm}\cdot\sum_i\,  s_{\dt_i} p_{r([v]_l,\dt_i)} s_{\dt_i}^*\cdot  
        s_{\gm}^*  \\ 
     & \leq \,   s_{\gm}  p_{[v]_l} s_\gm^*\sim p_{[v]_l}. 
\end{align*}   
Thus the projection $p_{[v]_l}$ is infinite 
since  $p_A$ is infinite by Proposition~\ref{QLE}.

On the other hand, 
every nonempty set in $\CE$ is a finite union of generalized 
vertices, hence it follows that 
every projection $p_B$, $B\in \CE$, is infinite. 
By Proposition~\ref{DHP}, 
any nonzero hereditary subalgebra of $C^*(E,\CL,\CE )$ contains 
a nonzero projection $p$ such that 
$s_\mu p_B s_\mu^*\preceq p$ for some $\mu\in \CL^*(E)$ 
and $B\in \CE $ with $B\subset r(\mu)$. 
But the projection $s_\mu p_B s_\mu^*\sim p_B$ is infinite, 
and thus $p$ is infinite.  

(b) We claim that if there is a vertex $[v]_l$ which does not 
connect to any loop, then $C^*(E,\CL,\CE )$ has an AF ideal.  

Suppose  a vertex $[v]_l \in \CE$ does not connect to any loop. 
By our assumption $(*)$, we may assume that 
$[v]_l$ is  minimal  in $\CE$. 
Let  
$$\{B_n : n\geq 1\}:=
\{B \in \CE_{\rm min} :B \subset r([v]_l,\af) 
\text{ for some } \af\in \CL^{\#}(E) \}.$$
(There could be only finitely many $B_n$'s, but 
the arguments given below work even in this case.) 
Note that each $B_n$ does not connect to any loop because 
$[v]_l$ would connect to a loop otherwise.
Hence for each $n\geq 1$ and $\af\in \CL(B_nE^{\geq 1})$, 
\begin{eqnarray}\label{Bn_noloop}
r(B_n,\af)\cap B_n=\emptyset
\end{eqnarray} 
since minimality of $B_n$ implies 
$B_n\subset r(B_n,\af)$ whenever $r(B_n,\af)\cap B_n\neq \emptyset$. 
If $n\neq m$ and 
$B_n\cap r(B_m, \ld) \neq \emptyset$ for some $\ld\in \CL^*(E)$ 
(hence, if  $B_n\subseteq r(B_m, \ld)$), 
we write $B_m>B_n$. 
Note that $B_n \not> B_n$ for all $n\geq 1$, 
and once $B_m>B_n$ holds for $n\neq m$, then
$B_n\not>B_m$ because of (\ref{Bn_noloop}).
The set of all finite unions of $B_n$'s which 
we denote by $H$ is readily seen 
to be a hereditary subset of $\CE$ by our assumption. 
With $p_n:=p_{B_n}$,  $n\geq 1$, the ideal $I_H$ 
generated by the projections $\{p_B:B\in H\}$ 
can be easily seen to be 
$$I_H=\overline{\rm span}\{ s_\af p_n s_\bt^* :  
\af,\bt\in \CL^*(E), \ n\geq 1\}.$$
Set 
$J_n:=\overline{\rm span}
\{s_\af p_n s_\bt^*: \af, \bt\in \CL^*(E)\}$ 
for each $n\geq 1$. 
Consider the product 
$X:=(s_\af p_n s_\bt^*)(s_\mu p_m s_\nu^*)$ 
of two elements 
$s_\af p_n s_\bt^*,$ and $s_\mu p_m s_\nu^*$ in $I_H$.
  Using the properties of the minimal sets $B_n$'s mentioned above, 
  one can check that 
\begin{eqnarray} \label{Jn-first}
\left\{           
  \begin{array}{ll}
   X\in  J_n, & \hbox{if  } n=m \\
   X\in  J_n,  & \hbox{if\ } B_m > B_n\\
   X\in  J_m,  & \hbox{if\ } B_n > B_m\\
   X=0,      & \hbox{otherwise.\ }
   \end{array}
   \right. 
\end{eqnarray} 
It then follows that both $J_n$ and 
$$I_n:=\overline{\rm span }\{ s_\af p_k s_\bt^*\, :
\, \af,\bt\in \CL^*(E),\ 1\leq k\leq n\}$$ 
are $C^*$-subalgebras of $I_H$, and moreover 
$I_H=\overline{\cup_{n\geq 1} I_n}$ holds. 
Thus it is enough to show that 
each $I_n$ is an AF algebra to obtain a contradiction. 

Fix $n\geq 1$ and consider the $C^*$-subalgebra $I_n$.
For each $k\in \{1,\dots, n\}$, 
$J_k$ is easily seen to be isomorphic to 
the $C^*$-algebra of compact operators 
$\CK(\ell^2(\CL(E^{\geq 1}B_k)))$ 
(see \cite[Corollary 2.2]{KPR}). 
Let $\{k_1,\dots, k_{n_k}\}\subset \{1,\dots, n\}$ 
be the set of all  such $k$'s that  
 $B_k\not>B_j$ for  all $1\leq j\leq n$. 
Then $J_k$ is an ideal of $I_n$ for each 
$k\in \{k_1,\dots, k_{n_k}\}$, and by (\ref{Jn-first}), 
$J_{k_i}J_{k_j}=\{0\}$ for $i\neq j$. 
Thus 
$I_n$ contains an AF ideal 
$\oplus_{k\in\{k_1,\dots,k_n\}} J_k$.
For the quotient map $\pi: I_n\to I_n/\oplus_{k} J_k$,  
the image $\pi(J_i)$ of $J_i$ is either zero or isomorphic to itself 
since each $J_i$ is a simple $C^*$-algebra.
Then one can apply the above argument to $\pi(I_n)$ 
which has the subalgebras $\pi(J_i) (\cong J_i)$,  
and find AF ideals and 
the quotient algebra that is generated by subalgebras 
each of which is isomorphic to the compact operators. 
This process should stop when a quotient algebra 
consists of ideals which are all isomorphic to 
the compact operators and 
mutually orthogonal, which 
proves that $I_n$ is an AF algebra. 
\end{proof}

\vskip 1pc

The following corollary 
generalizes \cite[Theorem 6.9]{BP2}. 
In fact, as mentioned earlier, 
the disagreeable and strongly cofinal labeled space 
$(E,\CL,\CE)$  in  Example~\ref{ex-VCL}  
satisfies the condition of Corollary~\ref{cor_loops} below,  
but it does not satisfy the condition of \cite[Theorem 6.9]{BP2}.

\vskip 1pc

\begin{cor}\label{cor_loops} 
Let $(E,\CL,\CE )$ be a strongly cofinal labeled space. 
If there exist a vertex $w \in E^0$, a strictly increasing sequence 
 $\{l_i\}_i$ of integers, and paths $\{\bt_i\}_i$
such that  $(\bt_i,[w]_{l_i})$ is a loop for each $i\geq 1$, 
then every vertex connects to 
a loop $(\bt_i,[w]_{l_i})$ for some $i\geq 1$. 

If, in addition, $(E,\CL,\CE )$ is disagreeable, then   
the simple $C^*$-algebra 
$C^*(E,\CL,\CE )$ is  purely infinite. 
In particular, if $(E,\CL,\CE )$ is a strongly cofinal and 
disagreeable labeled space with a cycle, then 
$C^*(E,\CL,\CE)$ is purely infinite and simple.
\end{cor}
 
\begin{proof}  
Given a generalized vertex $[v]_l$, 
since  $(E,\CL,\CE )$ is strongly cofinal,  
one can find $N\geq 1$ and paths $\dt_j$, $1\leq j\leq k$ 
such that 
$$[w]_{l_1}\subset r([w]_{l_1},\bt_1^N)\subset r(\bt_1^N) 
\subset \cup_{j=1}^k r([v]_l,\dt_j).$$
We may assume 
$w\in [w]_{l_1}\cap   r([v]_l,\dt_1)$, and then 
for all large $l_i$ we have 
$[w]_{l_i}\subset  r([v]_l,\dt_1)$ by (\ref{CE}). 
Thus $[v]_l$ connects to a loop $(\bt_{l_i},[w]_{l_i})$. 
If $(E,\CL,\CE )$ is disagreeable, by Theorem~\ref{thm-IH}.(a)
we see that the simple $C^*$-algebra 
$C^*(E,\CL,\CE)$ is purely infinite. 
\end{proof}
 
\vskip 1pc

\begin{remarks}\label{rmk of DHP} 
Let  $(E,\CL,\CE )$ be a labeled space. 
\begin{enumerate}
\item[(a)] One can modify the proof of Proposition~\ref{DHP}  to obtain the same 
criterion for quotient labeled spaces: 
Let $(E,\CL,\CE/R )$ be a  disagreeable quotient labeled space.
Then every nonzero hereditary $C^*$-subalgebra of $C^*(E,\CL,\CE/R)$ 
contains a nonzero projection $p$ such that 
$s_\mu p_{[A]} s_\mu^* \preceq p$ for some $\mu\in \CL_R^*(E)$ 
and $[A]\in \CE/R$ 
with $s_\mu p_{[A]} s_\mu^*\neq 0$ and $[A]\subset [r(\mu)]$. 

\item[(b)]
An argument similar to the proof of Theorem~\ref{thm-IH}.(a) 
shows the following: 
If a  quotient labeled space $(E,\CL,\CE/R)$ is disagreeable 
in which every vertex connects to a loop, 
then $C^*(E,\CL,\CE/R)$ has the property  (IH).
\end{enumerate}                    
\end{remarks}
 
\vskip 1pc

\subsection{An example: an infinite simple $C^*$-algebra 
of a labeled space $(E,\CL,\CE_\om)$ with no cycles}
It is well known that a simple graph $C^*$-algebra $C^*(E)$ is 
infinite if and if it is purely infinite, and this is 
also equivalent to the existence of a cycle (or a loop) in the graph $E$. 
Moreover, if this is the case, then every vertex in $E^0$ 
connects to a cycle.
 
For a simple labeled graph $C^*$-algebra $C^*(E,\CL,\CE)$ 
associated to a disagreeable and strongly cofinal labeled space 
$(E,\CL,\CE)$, we know that 
the existence of a cycle implies pure infiniteness of $C^*(E,\CL,\CE)$ 
by Corollary~\ref{cor_loops}.
But the converse may not be true. 
In fact,    
our purpose of this subsection is to provide an example of a  
labeled space $(E,\CL,\CE_\om)$ 
with no cycle but with a loop  
for which $C^*(E,\CL,\CE_\om)$ is purely infinite and simple. 
\vskip 1pc 

\begin{ex}\label{counterex} 
In the following  labeled graph $(E,\CL)$  

\vskip 1.5pc
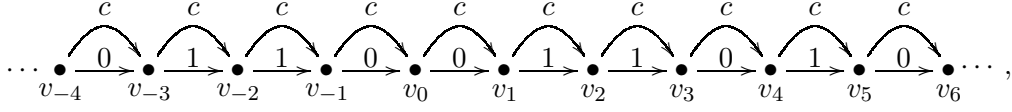
\xy  /r0.28pc/:(-44.2,0)*+{\cdots};(64.3,0)*+{\cdots,};
(-40,0)*+{\bullet}="V-4";
(-30,0)*+{\bullet}="V-3";
(-40,6)*+{ }="V-41";
(-20,0)*+{\bullet}="V-2";
(-10,0)*+{\bullet}="V-1"; (0,0)*+{\bullet}="V0";
(10,0)*+{\bullet}="V1"; (20,0)*+{\bullet}="V2";
(30,0)*+{\bullet}="V3";
(40,0)*+{\bullet}="V4";
(50,0)*+{\bullet}="V5";
(60,0)*+{\bullet}="V6";
(60,6)*+{  }="V61";
 "V-4";"V-3"**\crv{(-40,0)&(-30,0)}; ?>*\dir{>}\POS?(.5)*+!D{};
 "V-3";"V-2"**\crv{(-30,0)&(-20,0)}; ?>*\dir{>}\POS?(.5)*+!D{};
 "V-2";"V-1"**\crv{(-20,0)&(-10,0)}; ?>*\dir{>}\POS?(.5)*+!D{};
 "V-4";"V-3"**\crv{(-40,0)&(-35,10)&(-30,0)}; ?>*\dir{>}\POS?(.5)*+!D{};
 "V-3";"V-2"**\crv{(-30,0)&(-25,10)&(-20,0)}; ?>*\dir{>}\POS?(.5)*+!D{};
 "V-2";"V-1"**\crv{(-20,0)&(-15,10)&(-10,0)}; ?>*\dir{>}\POS?(.5)*+!D{};
 "V-1";"V0"**\crv{(-10,0)&(-5,10)&(0,0)}; ?>*\dir{>}\POS?(.5)*+!D{};
 "V0";"V1"**\crv{(0,0)&(5,10)&(10,0)}; ?>*\dir{>}\POS?(.5)*+!D{};
 "V1";"V2"**\crv{(10,0)&(15,10)&(20,0)}; ?>*\dir{>}\POS?(.5)*+!D{};
 "V2";"V3"**\crv{(20,0)&(25,10)&(30,0)}; ?>*\dir{>}\POS?(.5)*+!D{}; 
 "V3";"V4"**\crv{(30,0)&(35,10)&(40,0)}; ?>*\dir{>}\POS?(.5)*+!D{}; 
 "V4";"V5"**\crv{(40,0)&(45,10)&(50,0)}; ?>*\dir{>}\POS?(.5)*+!D{}; 
 "V5";"V6"**\crv{(50,0)&(55,10)&(60,0)}; ?>*\dir{>}\POS?(.5)*+!D{}; 
 "V-1";"V0"**\crv{(-10,0)&(0,0)}; ?>*\dir{>}\POS?(.5)*+!D{};
 "V0";"V1"**\crv{(0,0)&(10,0)}; ?>*\dir{>}\POS?(.5)*+!D{};
 "V1";"V2"**\crv{(10,0)&(20,0)}; ?>*\dir{>}\POS?(.5)*+!D{};
 "V2";"V3"**\crv{(20,0)&(30,0)}; ?>*\dir{>}\POS?(.5)*+!D{};
 "V3";"V4"**\crv{(30,0)&(40,0)}; ?>*\dir{>}\POS?(.5)*+!D{};
 "V4";"V5"**\crv{(40,0)&(50,0)}; ?>*\dir{>}\POS?(.5)*+!D{};
 "V5";"V6"**\crv{(50,0)&(60,0)}; ?>*\dir{>}\POS?(.5)*+!D{};
 (-35,1.5)*+{0};(-25,1.5)*+{1};
 (-15,1.5)*+{1};(-5,1.5)*+{0};(5,1.5)*+{0};
 (15,1.5)*+{1};(25,1.5)*+{1};(35,1.5)*+{0};
 (45,1.5)*+{1};(55,1.5)*+{0};
 (0.1,-2.5)*+{v_0};(10.1,-2.5)*+{v_1};
 (-9.9,-2.5)*+{v_{-1}};
 (-19.9,-2.5)*+{v_{-2}};
 (-29.9,-2.5)*+{v_{-3}};
 (-39.9,-2.5)*+{v_{-4}}; 
 (20.1,-2.5)*+{v_{2}};
 (30.1,-2.5)*+{v_{3}};
 (40.1,-2.5)*+{v_{4}}; 
 (50.1,-2.5)*+{v_{5}};
 (60.1,-2.5)*+{v_{6}};
(-35,7)*+{c};(-25,7)*+{c};(-15,7)*+{c};(-5,7)*+{c};
(5,7)*+{c};(15,7)*+{c};(25,7)*+{c};(35,7)*+{c};(45,7)*+{c};(55,7)*+{c}; 
\endxy 
\vskip 1.5pc
\noindent
the $\{0,1\}$-sequence is the bi-infinite Thue-Morse sequence  
which  is an infinite $\{0,1\}$-sequence with 
no finite subpaths of the form 
$\af\af\af_1$ for any
$\af=\af_1\cdots \af_{|\af|}\in \{0,1\}^{*}$  
(see \cite{GH}, \cite[Example 2.5, 2.7]{JKKP} for example).  
 Let $e_i\in E^1$, $i\in \mathbb Z$, be the edge of $E$    
with $s(e_i)=v_i$ and $\CL(e_i)\in \{0,1\}$. 
 For convenience, we write the Thue-Morse sequence as 
an infinite labeled path 
$$\om=\cdots \om_{-2}\om_{-1}\om_{0}\om_{1}\om_{2}\om_{3}\cdots,$$ 
where 
$\om_i=\CL(e_i)\in \{0,1\}$, $i\in \mathbb Z$. 
It is well known that $\om$ is not periodic but {\it almost periodic}
(or uniformly recurrent) in the sense that 
every finite subpath occurs with a bounded gap. 
We write the accommodating set $\CE$ by $\CE_\om$. 
Note that $(E,\CL,\CE_\om)$ is 
a disagreeable labeled space in which 
$(c^k,E^0)$ is a loop for each $k\geq 1$. 
Thus $C^*(E,\CL,\CE_\om)$ is an infinite $C^*$-algebra 
by  Proposition~\ref{QLE}, actually we have   
$$1=p_{r(c)}=s_c^*s_c\sim s_c s_c^*<s_cs_c^*+s_0s_0^*+s_1s_1^*=1.$$ 
We will show in the following proposition 
that the $C^*$-algebra $C^*(E,\CL,\CE_\om)$ is  simple,  
and 
the labeled space $(E,\CL,\CE_\om)$ contains no cycles 
and no  vertex $[v]_l\,(\neq E^0)$ connects to a loop. 
\end{ex} 

\vskip 1pc 

\begin{prop} Let $(E,\CL,\CE_\om)$  be the disagreeable labeled space 
of Example~\ref{counterex}. Then we have the following. 
\begin{enumerate}
\item[(a)] If $(\af,A)$ is a loop, then $\af=c^k$ for some $k\geq 1$.
\item[(b)] There exist no cycles.
\item[(c)] $C^*(E,\CL,\CE_\om)$ is purely infinite and simple.
\end{enumerate}
\end{prop}
\begin{proof}
(a) 
Suppose $(\af,A)$ is a loop with $k:=|\af|$. 
Then $\af^r$ occurs in $\CL^*(E)$ for every $r\geq 1$ 
since $A\subset r(A,\af)\subset r(A,\af^2)\subset \cdots$. 
If $\af$ has a $\{0,1\}$-subpath, say $\af'$ (which 
we may assume ends with $0$), then 
$r(\af')$ must contain a periodic finite set of arbitrary cardinality. 
But this is not possible by Gelfond's result \cite{G} saying that 
for any $a,b\in \mathbb N$, 
$$\Big |\big\{ n<N: \om_{an+b}=0\big\}\Big|=\frac{N}{2}+O(N^\ld).$$
Thus we have $\af=c^k$. 
 
(b) 
Suppose $(\af,A)$ is a cycle in $(E,\CL,\CE_\om)$.
Then $\af=c^k$ for some $k\geq 1$ by (a), and  
$B= r(B,c^k)$ holds for any nonempty subset $B\subset A$.
Especially for $B=r(\bt)\subset A$ where $\bt$ is a 
$\{0,1\}$-path  of length $|\bt|>3k$, 
we have 
$$r(\bt)=B= r(B,c^k)=r(r(\bt),c^k)=r(\bt c^k)
=\underset{\dt\in \{0,1\}^k}{\cup}\,r(\bt\dt).$$ 
Choose any vertex $v\in r(\bt)$. Then 
$v\in r(\bt\dt)$ for some $\{0,1\}$-path $\dt$ with $|\dt|=r$. 
This means that there must be a path $\bt'$ such that 
$\bt'\bt=\bt\dt$. 
Since $|\bt|>3|\dt|$, we obtain that 
$\bt$ ends with $\dt^3$, which is not possible 
since  $\bt$ is a subpath of the Thue-Morse sequence $\om$. 

(c) 
We first show that $C^*(E,\CL,\CE_\om)$ has the 
property (IH). 
Let $[v]_l$  be a generalized vertex. 
Then there is a finite $\{0,1\}$-path $\bt\in \CL(E^l)$ such that 
$[v]_l=r(\bt)$.
Since $\bt$ is a subpath of the almost periodic sequence 
$\om$,  it occurs with a bounded gap, say $L_\bt$. 
Thus we can find finitely many paths 
$\dt_i$ of the form $\dt_i:=\af_i c$ for 
some $\af_i\in \{0,1\}^*$ such that 
$s_{\dt_i}^* s_{\dt_j}=0$ for $i\neq j$ and 
$$E^0\subset \cup_i\, r([v]_l,\dt_i)=\cup_i \, r(\bt\dt_i). $$
This is possible since if  $w\in E^0$ is a vertex, then 
$w\in r(\bt\af c)$ for some $\{0,1\}$-path 
$\af$ (possibly $\af=\epsilon$) 
such that $|\af|\leq L_\bt$.
Thus $[v]_l$ connects to a loop $(c,E^0)$, 
and by Theorem~\ref{thm-IH}.(a)
$C^*(E,\CL,\CE_\om)$ has the property (IH). 

To see that $C^*(E,\CL,\CE_\om)$ is simple,  
we only need to show by Remark~\ref{simple condition} 
that $(E,\CL,\CE_\om)$ is strongly cofinal. 
Let $[v]_l=r(\bt)$, $\bt\in\{0,1\}^l$, be a generalized vertex 
as above and  let
 $x\in \overline{\CL(E^\infty)}$. 
We already  know that 
$E^0\subset \cup_i\, r([v]_l,\dt_i)$ for finitely many paths 
$\dt_i$'s, and thus  
$\displaystyle r(x_{[1,N]})\subset E^0
=  \cup_{ i} r([v]_l,\dt_i)$ 
for any  $N\geq 1$, which  
shows that $(E,\CL,\CE_\om)$ is strongly cofinal as desired.
\end{proof}

\vskip 1pc 

\begin{remark} 
Pure infiniteness of $C^*(E,\CL,\CE_\om)$ can also be obtained if 
one can prove that the standard action of the inverse semigroup  
$$\CS:=\big\{ s_\af p_{r(\dt)}s_\bt^*: 
\af,\bt,\dt\in \CL^*(E),\,  |\dt|\geq \max\{|\af|,|\bt|\}\big\}$$
associated to the labeled space $(E,\CL,\CE_\om)$ 
on the tight spectrum $\hat\CE(\CS)_{tight}$ of $\CE(\CS)$   
is locally contracting in the sense of 
\cite[Definition 6.2]{EP} (although this does not seem 
to work easily), because then by 
\cite[Proposition 6.3]{EP} 
the tight groupoid $\CG_{tight}(\CS)$  is 
locally contracting, and hence by 
\cite[Proposition 2.4]{A} 
one sees that 
the $C^*$-algebra $C^*(E,\CL,\CE_\om)$ is purely infinite. 

For this inverse semigroup $\CS$ of $C^*(E,\CL,\CE_\om)$, 
 the standard action is locally 
contracting  if and only if 
$\CS$ is locally contracting in the sense of 
\cite[Definition 6.4]{EP} 
by \cite[Theorem 6.5]{EP} since 
every tight-filter is a ultra-filter which we briefly show here.  
First recall that a filter $\xi$ is a {\it ultra-filter} if and only if 
$\xi$ contains every idempotent $f$ such that 
$fe\neq 0$ for every $e\in \xi$. 
Let $\eta$ be a tight-filter and 
$f=s_\sm p_{r(\bt)} s_\sm^*\in \CE(\CS)$ be an idempotent 
such that $fe\neq 0$ for every $e\in \eta$. 
There is a unique  tight character $\phi$ on $\CE(\CS)$ such that 
$\eta=\eta_\phi:=\{e\in \CE(\CS): \phi(e)=1\}$. 
Pick $e=s_\mu p_{r(\af)} s_\mu^*\in \eta_\phi$.  
By \cite[Proposition 1.8]{E} we have for every finite cover 
$Z\subset \CE(\CS)$ for $e$, 
$$ \underset{z\in Z}{\vee} \phi(z)\geq \phi(e)=1,$$
which implies that if 
$e=\sum_{i=1}^n e_i$ is a finite sum of mutually orthogonal 
idempotents $e_i\in \CE(\CS)$,  
for example, one can write $e$ as follows 
$$e= \sum_{\nu\in \CL^k(E)} s_{\mu\nu}p_{r(\af\nu)} s_{\mu\nu}^*   
\ \text{ or }\
e=\sum_{\af_i\in \{0,1\}^l}  s_\mu p_{r(\af_i\af)} s_\mu^*$$ 
for $k,l\geq 1$, 
then 
$\phi(e_i)=1$ for one and only one $e_i$ because 
$\{e_1, \dots, e_n\}$ is a cover for $e$ and 
$\phi$ is a multiplicative function with $\phi(0)=0$.  
Also the identity 
$r(\af)=
\underset{{\substack{\af_i\in \{0,1\}^{|\af|}\\ r(\af_i)=r(\af)}}}{\cup} r(\af_i)$ 
allows us to assume that $\af$ is a $\{0,1\}$-path. 
Thus we may pick $e\in \eta_\phi$ such that 
$|\mu|\geq |\sm|$ and $\af$ is a  $\{0,1\}$-path
with $|\af|\geq |\bt|+|\mu|-|\sm|$. 
Then 
\begin{align*} 
0\neq fe &= s_\sm p_{r(\bt)} s_\sm^*\cdot s_\mu p_{r(\af)}s_\mu^*\\
&=\sum_{|\sm'|=|\mu|-|\sm|} s_{\sm\sm'} p_{r(\bt\sm')} s_{\sm\sm'}^*\cdot 
s_\mu p_{r(\af)}s_\mu^* \\
&= s_{\sm\sm'} \, p_{r(\bt\sm')}\cdot p_{r(\af)} s_\mu^*\ (\text{for }\sm\sm'=\mu)\\
&= s_{\sm\sm'}\cdot\, \sum_{\bt_i\bt \sm_j\in \{0,1\}^{|\af|}}
 p_{r(\bt_i\bt \sm_j)}\cdot p_{r(\af)} s_\mu^* \\
&=  s_{\sm\sm'}p_{r(\af)}  s_\mu^* \\
& =e\in \eta_\phi.
\end{align*}
Since $f\geq fe$ and $\eta_\phi$ is a filter, we have 
$f\in \eta_\phi$, which proves that 
$\eta_\phi$ is a ultra-filter.
\end{remark}

\vskip 1pc

\section{Purely infinite labeled graph $C^*$-algebras}

In this section, we find equivalent conditions for 
a labeled graph $C^*$-algebra $C^*(E,\CL,\CE)$ to be purely infinite 
when the labeled space $(E,\CL,\CE)$ is {\it strongly disagreeable} 
(see Definition~\ref{ddef-st disagreeable}).

We begin with the following Cuntz-Krieger uniqueness theorem 
for quotient labeled graph $C^*$-algebras. 
Although the theorem can be obtained 
by standard arguments as the proof of  \cite[Theorem 5.5]{BP2} 
for labeled graph $C^*$-algebras, 
we give a proof here
only to make sure that  our definition of 
a `disagreeable' quotient labeled space (Definition~\ref{disagreeable q}) 
does not cause any degenerate relations. 
We will use this uniqueness theorem later in Lemma~\ref{GI}.
 
\vskip 1pc

\begin{thm}\label{QUT}
Let  $(E,\CL,\CE/R)$ be a disagreeable quotient labeled space  
and let $C^*(E,\CL,\CE/R)=C^*(s_a,p_{[A]})$. 
If $\{S_a,P_{[A]}\}$ is a representation of $(E,\CL,\CE/R)$ such that 
$P_{[A]} \neq 0$ for $[A] \neq [\emptyset]$ and 
$S_a \neq 0$ for $[r(a)] \neq [\emptyset]$, 
then the homomorphism 
$$\pi_{S,P}:C^*(E,\CL,\CE/R) \rightarrow C^*(S_a,P_{[A]})$$ 
with $\pi_{S,P}(s_a)=S_a$ and 
$\pi_{S,P}(p_{[A]})=P_{[A]}$ is faithful.
\end{thm}
 
\begin{proof} 
As usual (for example, see  \cite[Theorem 5.5]{BP2}), 
it is enough to show that 
\begin{enumerate}
\item[(a)] $\pi_{S,P}$ is faithful on $C^*(E,\CL,\CE/R)^{\gm}$ and
\item[(b)] $\|  \pi_{S,P}(\Phi(a))\| \leq \|
\pi_{S,P}(a)\|$  for all $a \in C^*(E,\CL, \CE/R)$,  
where $\Phi$ is the conditional expectation onto the fixed point algebra 
$C^*(E,\CL,\CE/R)^{\gm}$.
\end{enumerate} 
But (a) was shown in the proof of \cite[Theorem 4.2]{JKP}, and 
thus we show (b). 
By (\ref{q-span}),  it suffices to
prove (b) for an element  $a$ of the form
$$ a = \sum_{(\af, [[w]_l], \bt) \in F} \, 
c_{\af, [[w]_l], \bt}\,s_{\af}p_{[[w]_l]}s_{\bt}^*, $$
where $F$ is a finite set and $c_{\af, [[w]_l], \bt}\in \mathbb C$. 
Let $k = \max \{|\af|,|\bt| : (\af, [[w]_l], \bt) \in F  \}$. 
We can choose $k \in \mathbb{N} $ by applying Definition~\ref{QR}(iv)  
and changing $F$ if necessary 
such that $\min \{|\af|,|\bt| \} = k$ for
every $ (\af, [[w]_l], \bt) \in F$.  Let  $M = \max\{|\af|,|\bt| :
(\af, [[w]_l], \bt) \in F  \}$. 
We also may suppose that $l \geq M-k$. 
 Since $|\af|=|\bt|=k$ whenever $|\af|=|\bt|$, we have 
 $$ \Phi(a) = \sum_{\substack{(\af, [[w]_l], \bt) \in F\\ |\af|=|\bt|=k}}
 c_{\af, [[w]_l], \bt}\, s_{\af}p_{[[w]_l]}
 s_{\bt}^* \in \oplus_{[[w]_l]} \CF^k([[w]_l]).$$ 
 Thus the norm 
 $\|\pi_{S,P}(\Phi(a))\|$ is attained on a direct summand, 
say $\CF^k([[v]_l])$. 
For this $[[v]_l]$, let 
$F_{[v]_l}:=\{(\af, [[v]_l], \bt): (\af, [[v]_l], \bt) \in F\}$. 
 Then we have 
 $$\|\pi_{S,P}(\Phi(a))\|
 =\Big\| \sum_{\substack{(\af, [[v]_l], \bt) \in F_{[v]_l}\\ |\af|=|\bt|=k}}
 c_{\af, [[v]_l], \bt}\, S_{\af}P_{[[v]_l]}S_{\bt}^*\, \Big\|.$$
It is easy to see  that the $*$-algebra
${\rm span} \{S_{\af}P_{ [[v]_l]}S_{\bt}^* : \af, \bt \in G\}$, where  
$$G:=\{\af \in \CL_R(E^k) : \text{either } (\af, [[v]_l], \bt) \in F_{[v]_l} 
 \text{ or } (\bt, [[v]_l], \af) \in F_{[v]_l} \text{ with } 
|\af|=|\bt|\,\},$$  
is a finite dimensional matrix algebra and 
contains the element 
$$b_v:=\sum_{\substack{(\af, [[v]_l], \bt) \in F_{[v]_l}\\ |\af|=|\bt|=k}}
c_{\af, [[v]_l], \bt}\, S_{\af}P_{[[v]_l]}S_{\bt}^*.$$
Since $[[v]_l]$ is disagreeable (see Definition~\ref{disagreeable q}),
we can choose a path 
 $\ld \in \CL_R([[v]_l]E^{>M})$ 
which has no
factorization $\ld = \ld' \ld'' =\ld''\gm$ for 
$\ld',\gm \in \CL_R^{*}(E)$ with $|\ld'|=|\gm| \leq M-k \leq l$  
and $[r([v]_l,\ld)] \neq [\emptyset]$. 
This gives a nonzero projection 
 $$ Q = \sum_{\nu \in G}S_{\nu\ld}P_{[r([v]_l,
\ld)]}S_{\nu\ld}^*.$$
It is routine to check the following 
(for example, see \cite[Theorem 5.5]{BP2})
\begin{enumerate}
\item[(i)]  $\| Q \,\pi_{S,P}(\Phi(a))\,Q \| =\|\pi_{S,P}(\Phi(a)) \| $
\item[(ii)] $Q\,S_{\af}P_{[[v]_l]}S_{\bt}^*\,Q = 0$ when 
$(\af, [[v]_l],\bt) \in F$  and $|\af| \neq |\bt|$.
\end{enumerate}
Then we see that \begin{align*}  
\|  \pi_{S,P}(\Phi(a))\| &=\|Q  \pi_{S,P}(\Phi(a) ) Q\,\Big\| \\
& =\Big\|\, Q\Big(\sum_{\substack{(\af, [[v]_l], \bt) \in F_{[v]_l}\\ 
|\af|=|\bt|}} c_{\af, [[v]_l], \bt}\,S_{\af}P_{[[v]_l]}S_{\bt}^*\Big)Q\,\Big\|\,  \\
& =\Big\|Q\Big(\sum_{(\af, [[v]_l], \bt) \in F_{[v]_l}} 
   c_{\af, [[v]_l], \bt}\,S_{\af}P_{[[v]_l]}S_{\bt}^*\Big)Q\| \\
& \leq \Big\|\sum_{(\af, [[v]_l], \bt) \in F_{[v]_l}}
   c_{\af, [[v]_l], \bt}\,S_{\af}P_{[[v]_l]}S_{\bt}^*\,\Big\|\\ 
& = \|\pi_{S,P}(a)\|.
\end{align*} 
This completes the proof.
\end{proof}

\vskip 1pc

\subsection{Gauge-invariant 
Ideals of $C^*$-algebras of strongly disagreeable labeled spaces}
Here we introduce 
a property, called `strongly disagreeable', 
 of a labeled space $(E,\CL,\CE)$ 
which we believe is a necessary condition
for the $C^*$-algebra $C^*(E,\CL,\CE)$ to be  purely infinite,  
at least in many practical examples we easily come across with 
including all the graph $C^*$-algebras.

As is seen from the following labeled graph, 
a $C^*$-algebra $C^*(E,\CL,\CE)$ is 
not necessarily purely infinite 
for a disagreeable labeled space  $(E,\CL,\CE)$ 
if it has a non-disagreeable quotient labeled space: 

\vskip 1pc \hskip 7pc 
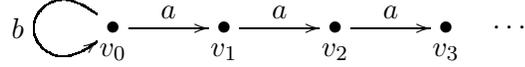
\xy /r0.5pc/:
 (20,0)*+{\bullet}="V1";
 (27,0)*+{\bullet}="V2";
 (34,0)*+{\bullet}="V3";
 (41,0)*+{\bullet}="V4";
 "V1";"V1"**\crv{(20,0)&(17,3)&(14,0)&(17,-3)&(20,0)};?>*\dir{>}\POS?(.5)*+!D{};
 "V1";"V2"**\crv{(20,0)&(27,0)};?>*\dir{>}\POS?(.5)*+!D{};
 "V2";"V3"**\crv{(27,0)&(34,0)};?>*\dir{>}\POS?(.5)*+!D{};
 "V3";"V4"**\crv{(34,0)&(41,0)};?>*\dir{>}\POS?(.5)*+!D{};
(14,0)*+{b}, (23.5,1)*+{a},(30.5,1)*+{a},(37.5,1)*+{a}, 
(20,-1.5)*+{v_0},(27,-1.5)*+{v_1},(34,-1.5)*+{v_2},(41,-1.5)*+{v_3}, 
(45,0)*+{\cdots},
\endxy

\vskip 2pc   

\begin{dfn} \label{ddef-st disagreeable}
We say that $(E,\CL,\CE)$ is {\it strongly disagreeable} if 
the quotient labeled space $(E,\CL, \CE/H)$ is disagreeable  
for every hereditary saturated subset $H$ of $\CE$.
\end{dfn}

\vskip 1pc  

\begin{remarks} Let $E$ be a row-finite graph with no sinks.
\begin{enumerate}
\item[(a)] $E$  satisfies condition (K) if and only if 
the labeled space $(E,\CL_{id}, \CE)$ is strongly disagreeable: 
Recall first that 
$E$ satisfies condition (K) if and only if 
every ideal of $C^*(E)$ is gauge-invariant 
(see \cite[Lemma 2.2]{BHRS} for the fully general version). 
Thus Proposition~\ref{GI} below implies that 
if $(E,\CL_{id}, \CE)$ is strongly disagreeable, 
then $E$ has condition (K).
The converse can be easily seen by considering 
a quotient graph because 
a hereditary saturated set $H$ ($\CE/H$, respectively) 
of $\CE$ can be identified with a  
hereditary saturated set of $E^0$ ($E^0\setminus H$, respectively). 

\item[(b)] A {\it strongly aperiodic}  higher rank graph was   
introduced in \cite{KP}, and 
it was shown \cite[Corollary 3.9]{PSS} that 
a $k$-graph $\Lambda$ is strongly aperiodic if and only if 
every ideal of the $C^*$-algebra $C^*(\Lambda)$ associated to $\Ld$ is 
gauge invariant. 
Thus strong aperiodicity can be viewed 
as the higher-rank graph analogue for condition  (K).
\end{enumerate}
\end{remarks}
 
\vskip 1pc

\begin{prop}\label{GI}  
If $(E,\CL,\CE )$ is strongly disagreeable,  
every ideal  of $C^*(E,\CL,\CE)$ is gauge-invariant.
\end{prop}
 
\begin{proof}
Let $I$ be an ideal of $C^*(E, \CL,\CE )$. 
Then  $H_I=\{A \in \CE :p_A \in I\}$ is a 
saturated hereditary subset of $\CE $  and 
the ideal $I_{H_I}$ generated by the projections 
$\{p_A : A \in H_I\}$ is gauge invariant 
\cite[Lemma 3.9]{JKP}. 
Since $I_{H_I} \subseteq I$, the quotient map 
$$q: C^*(E,\CL,\CE)/I_{H_I} \rightarrow C^*(E,\CL,\CE )/I$$ 
given by $q(s+I_{H_I}):=s+I$ for $s\in C^*(E,\CL,\CE)$, 
is well-defined.
From \cite[Theorem 5.2]{JKP}, we have an isomorphism 
$\pi: C^*(E,\CL,\CE/H_I)\to C^*(E,\CL,\CE )/I_{H_I}$  
which maps the canonical generators to  the canonical generators.   
Then the composition map  
$q \circ \pi :C^*(E,\CL,\CE/H_I) \rightarrow  C^*(E,\CL,\CE )/I$ 
satisfies  
\begin{align*}
q \circ \pi(p_{[A]})& =q(p_A+I_{H_I})=p_A+I\\
q \circ \pi(s_a) & =q(s_a+I_{H_I})=s_a+I 
\end{align*} 
for $[A]\in \CE/H_I$ and $a\in \CA_{H_I}$. 
If $p_{[A]}\neq 0$, then  $[A]\neq [\emptyset]$ in $\CE/H_I$, hence 
 $A \notin H_I$.
Thus  $p_A+I\in C^*(E,\CL,\CE)/I$ is a nonzero projection.
 Also for $a\in \CA_{H_I}$ 
(namely $[r(a)]\neq [\emptyset]$ in $H_I$), we have 
  $q \circ \pi(s_a)=q(s_a+I_{H_I})=s_a+I \neq 0$ 
 because otherwise $p_{r(a)}=s_a^* s_a\in I$, a contradiction to  
 $[r(a)] \neq [\emptyset]$ in $H_I$.
Since the quotient labeled space $(E, \CL, \CE/H_I)$ is disagreeable 
by our assumption, 
we see that the map $q \circ \pi $ is injective by  Theorem~\ref{QUT}.  
   Thus  $q$ is be injective, so that $I$ must coincide with 
   the gauge invariant ideal  $I_{H_I}$. 
       
   Note that if $I$ is an ideal such that $H_I=\{ \emptyset \}$, 
   then  $I_{H_I}=\{0\}$, and   
   $q \circ \pi$ is the quotient map $q: C^*(E,\CL,\CE) \rightarrow C^*(E,\CL,\CE)/I$. 
   Since the family $\{p_A+I,\, s_a+I \}$ is a 
   representation of the labeled space $(E,\CL,\CE)$ 
   in the $C^*$-algebra $C^*(E,\CL,\CE)/I$ such that 
   $p_A+ I \neq 0$ and $s_a+I \neq 0$ for each $A \in \CE $ and $a \in \CA$,  
  the disagreeability of $(E,\CL, \CE)$ asserts 
by the Cuntz-Krieger unique theorem that $q$ is injective.  
  Thus we have $I=\{0\}$.  
\end{proof}

\vskip 1pc

\begin{ex} For the following labeled graph 
\vskip 1pc
\hskip 7pc 
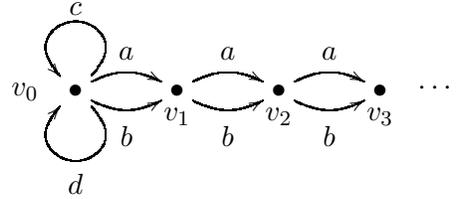
\xy /r0.4pc/:
(-24,0)*+{\bullet}="V-3"; 
(-16,0)*+{\bullet}="V-2";
(-8,0)*+{\bullet}="V-1"; 
(0,0)*+{\bullet}="V0";
"V-3";"V-3"**\crv{(-24,0)&(-20,3.3)&(-24,6.6)&(-28,3.3)&(-24,0)};?>*\dir{>}\POS?(.5)*+!D{};
"V-3";"V-3"**\crv{(-24,0)&(-20,-3.3)&(-24,-6.6)&(-28,-3.3)&(-24,0)};?>*\dir{>}\POS?(.5)*+!D{};
 "V-3";"V-2"**\crv{(-24,0)&(-20,3)&(-16,0)};
 ?>*\dir{>}\POS?(.5)*+!D{};
 "V-3";"V-2"**\crv{(-24,0)&(-20,-3)&(-16,0)};
 ?>*\dir{>}\POS?(.5)*+!D{}; 
 "V-2";"V-1"**\crv{(-16,0)&(-12,3)&(-8,0)};
 ?>*\dir{>}\POS?(.5)*+!D{};
 "V-2";"V-1"**\crv{(-16,0)&(-12,-3)&(-8,0)};
 ?>*\dir{>}\POS?(.5)*+!D{};
 "V-1";"V0"**\crv{(-8,0)&(-4,3)&(0,0)};
 ?>*\dir{>}\POS?(.5)*+!D{};
 "V-1";"V0"**\crv{(-8,0)&(-4,-3)&(0,0)};
 ?>*\dir{>}\POS?(.5)*+!D{};
(5,0)*+{\cdots ,};
 (-20,3)*+{a};(-12,3)*+{a};(-4,3)*+{a};
 (-20,-3.5)*+{b};(-12,-3.5)*+{b};(-4,-3.5)*+{b};
 (-24,6.5)*+{c};(-24,-7.3)*+{d};
 (-28,0)*+{v_0};(-16,-2)*+{v_1};(-8,-2)*+{v_2};(0,-2)*+{v_3};
\endxy

\vskip 1pc

\noindent 
it is not hard to see that 
$\CE =\{A \subset E^0 : A  \text{ is finite  or co-finite}\,\}$.   
The labeled space $(E,\CL,\CE)$ is disagreeable, and 
 there are only three non-trivial hereditary saturated subsets of 
 $\CE $; $H_1=\{A \subset E^0  : |A |<\infty \}$, 
 $H_2=\{A \subset E^0\setminus \{v_0\} : |A |<\infty \text{ or } |A^c |<\infty \}$,  
 and $H_3=\{A \subset E^0\setminus \{v_0\} : |A |<\infty\}$. 
Note that $H_3\subset H_2\cap H_1$. 
To see what each quotient algebra 
$C^*(E,\CL,\CE/H_i)$ would be like, 
it is helpful to view each quotient labeled space $(E,\CL, \CE/H_i)$ 
as if it is a labeled space:
\vskip 1pc
\hskip 3.5pc
\xy /r0.38pc/:
(0,0)*+{\bullet}="u";
  "u";"u"**\crv{(0,0)&(4,3.5)&(8,0)&(4,-3.5)&(0,0)};?>*\dir{>}\POS?(.5)*+!D{};
  "u";"u"**\crv{(0,0)&(-4,3.5)&(-8,0)&(-4,-3.5)&(0,0)};?>*\dir{>}\POS?(.5)*+!D{};
   (-8,0)*+{a},(8,0)*+{b},(-20,0)*+{(E,\CL, \CE/H_1):},(0,-3.8)*+{[E^0]}
   \endxy
\vskip 1pc
\hskip 3.5pc
\xy /r0.38pc/:
(0,0)*+{\bullet}="u";
  "u";"u"**\crv{(0,0)&(4,3.5)&(8,0)&(4,-3.5)&(0,0)};?>*\dir{>}\POS?(.5)*+!D{};
   "u";"u"**\crv{(0,0)&(-4,3.5)&(-8,0)&(-4,-3.5)&(0,0)};?>*\dir{>}\POS?(.5)*+!D{};
   (-8,0)*+{c},(8,0)*+{d},(-20,0)*+{(E,\CL,\CE/H_2):},(0,-3.8)*+{[\{v_0\}]}
   \endxy
\vskip 1pc
\hskip 3.5pc
\xy /r0.38pc/:
(0,0)*+{\bullet}="u";
(0,0)*+{\bullet}="u";
  "u";"u"**\crv{(0,0)&(4,3.5)&(8,0)&(4,-3.5)&(0,0)};?>*\dir{>}\POS?(.5)*+!D{};
   "u";"u"**\crv{(0,0)&(-4,3.5)&(-8,0)&(-4,-3.5)&(0,0)};?>*\dir{>}\POS?(.5)*+!D{};
   (-8,0)*+{c},(8,0)*+{d},(-20,0)*+{(E,\CL,\CE/H_3):},(0,-3.8)*+{[\{v_0\}]};  
(20,0)*+{\bullet}="v";
  "v";"v"**\crv{(20,0)&(16,3.5)&(12,0)&(16,-3.5)&(20,0)};?>*\dir{>}\POS?(.5)*+!D{};
   "v";"v"**\crv{(20,0)&(24,3.5)&(28,0)&(24,-3.5)&(20,0)};?>*\dir{>}\POS?(.5)*+!D{};
   (12,0)*+{a},(28,0)*+{b},(20,-4)*+{[r(a)]};   
     \endxy
\vskip 1pc
\noindent
Then it is rather obvious that $(E,\CL,\CE/H_i)$ is disagreeable 
for all $i=1,2,3$, and  thus
  by Proposition~\ref{GI} 
every ideal of $C^*(E,\CL,\CE )$ is gauge-invariant. 
Consequently, the $C^*$-algebra $C^*(E,\CL,\CE )$ has only 
three non-trivial ideals $I_{H_i}$, $i=1,2,3$, by 
Remark~\ref{remark-review}.(a).
\end{ex}
\vskip 1pc
 
\subsection{Purely infinite  $C^*$-algebras of strongly disagreeable 
labeled spaces} 

We see from the following lemma 
that if $C^*(E,\CL,\CE)$ is purely infinite, 
every loop on a minimal set  should have an exit. 
 
\vskip 1pc

\begin{lem}\label{MLQ} 
Let $(E,\CL,\CE )$ be a labeled space 
with a minimal set $A\in \CE_{\rm min}$. 
If there is a loop $(\af,A)$ which has no exits, then 
the $C^*$-algebra $C^*(E,\CL,\CE )$ 
contains a hereditary subalgebra 
isomorphic to $M_n(C(\mathbb{T}))$ 
for some $n \in \mathbb{N}$.
\end{lem}

\begin{proof} 
Since $(\af, A)$ has no exits,   
we have $A=r(A,\af)$ and $\CL(AE^{|\af|})=\{\af\}$. 
These two conditions imply that   
\begin{eqnarray}\label{ne}  
\CL(AE^{\geq 1}) =   \{ \alpha ^k\af' : k \geq 0, \af' 
~\text{is an initial path of}~ \af \}. \end{eqnarray}  
Thus 
$$p_A=s_{\af_1}p_{r(A,\af_1)}s_{\af_1}^* = 
s_{\af_{[1,2]}}p_{r(A,\af_{[1,2]})}s_{\af_{[1,2]}}^*
=\cdots=s_{\af}p_{r(A,\af)}s_{\af}^*.$$
Furthermore the relative ranges $r(A,\af_{[1,i]})$  
are all minimal and mutually disjoint for $1\leq i\leq |\af|$. 
In fact, 
if $r(A,\af_{[1,i]})$ is the union of two disjoint sets 
$B_1, B_2 \in \CE $, 
then
$$A=r(A,\af)=r(r(A,\af_{[1,i]}),\af_{[i+1,|\af|]})
=r(B_1,\af_{[i+1,|\af|]} ) \cup r(B_2,\af_{[i+1,|\af|]} ),$$
and $r(B_1,\af_{[i+1,|\af|]} ) \cap r(B_2,\af_{[i+1,|\af|]} ) 
= \emptyset $ since 
$(E,\CL,\CE )$ is weakly left-resolving. 
From minimality of $A$, one of the two sets must be empty. 
But this also implies that $B_1= \emptyset$ or $B_2= \emptyset$ because 
$r(A,\af_{[1,i]})=B_1 \cup B_2$ has no sinks. 
Hence, $r(A,\af_{[1,i]})$ is minimal. 
We may assume that $\af$ is a loop at $A$ with the smallest length 
(that is, $|\af| \leq |\bt|$ for any loop 
$(\bt,A)$ or $\bt \in \CL(AE^{\geq 1}A)$), 
and then  we can easily see that 
$$r(A,\af_{[1,i]}) \cap r(A,\af_{[1,j]}) = \emptyset$$
for $1 \leq i,j \leq |\af|$, $i \neq j$.
For $1 \leq i \leq n :=|\af|$, set $A_i:=r(A,\af_{[1,i]})$  and 
$$p:= p_{A_1}+\cdots+p_{A_n}.$$ 
Now we claim that $pC^*(E,\CL,\CE )p\cong C(\mathbb{T}) \otimes M_n$. 
For $s_{\mu}p_Cs_{\nu}^* \in C^*(E,\CL,\CE )$, if
$$p(s_{\mu}p_Cs_{\nu}^*)p = 
\sum_{i,j}s_{\mu}p_{r(A_i,\mu)\cap C \cap r(A_j,\nu)}s_{\nu}^* \neq 0,$$ 
then
$s_{\mu}p_{r(A_i,\mu)\cap C \cap r(A_j,\nu)}s_{\nu}^* \neq 0$ for some $i,j$. 
Hence 
$$r(A_i,\mu)\cap  r(A_j,\nu) 
=r(A,\af_{[1,i]}\mu)\cap r(A,\af_{[1,j]}\nu)\neq \emptyset.$$ 
We then see from (\ref{ne}) that 
the paths $\mu$, $\nu$ are of the form 
$$\mu=\af_{[i+1, ,n]}\af^l\af_{[1,k]},\ 
\nu=\af_{[j+1,n]}\af^m\af_{[1,k']}$$ 
for some $l,m \geq 0$ and $0\leq k, k'<n$. 
Then $\emptyset\neq r(A,\af_{[1,i]}\mu)\cap r(A,\af_{[1,j]}\nu) 
= r(A,\af^{l+1}\af_{[1,k]})\cap r(A,\af^{m+1}\af_{[1,k']})
=  r(A, \af_{[1,k]})\cap r(A, \af_{[1,k']})
=  A_k\cap A_{k'}$. 
Thus,  $k=k'$ (namely, $r(A_i,\mu) =r(A_j,\nu)$)  and 
\begin{align} 
s_{\mu}p_{r(A_i,\mu)\cap C \cap r(A_j,\nu)}s_{\nu}^* 
& =s_{\af_{[i+1,n]}\af^l\af_{[1,k]}}
  p_{r(A,\af_{[i,k]})}s_{\af_{[j+1,n]}\af^m\af_{[1,k]}}^* \nonumber \\
& = s_{\af_{[i+1,n]}\af^l\af_{[1,k]}}p_{A_k}s_{\af_{[j+1,n]}\af^m\af_{[1,k]}}^* . 
\nonumber \\ \nonumber 
\end{align}  
This means that the hereditary subalgebra  
$pC^*(E,\CL,\CE )p$ is generated by the elements 
$p_{A_{i-1}}s_{\af_i}(=s_{\af_i}p_{A_i})$. 
Let $\gm$ be the restriction of the gauge action on 
$C^*(E,\CL,\CE )$ to the hereditary algebra 
$pC^*(E,\CL,\CE )p$ (which is obviously gauge-invariant) 
and let $\bt$ be the gauge action of the universal $C^*$-algebra 
$C(\mathbb{T}) \otimes M_n$ (which is actually a graph $C^*$-algebra) 
generated by 
the partial isometries $t_1,\dots,t_n$ 
satisfying
$$t_i^*t_i=t_{i+1}t_{i+1}^*,\, t_n^*t_n=t_1t_1^*,\, 
\text{ and }\, \sum_{j=1}^{n}t_j^*t_j=1$$
for $1 \leq i \leq n-1$. 
But the elements $s_{\af_i} p_{A_i}$, $1\leq i\leq n$, 
satisfy the relations with $p$ in place of $1$, hence 
by universal property of $C(\mathbb{T})\otimes M_n$, 
there exists a  homomorphism
  $$\pi:C(\mathbb{T}) \otimes M_n \rightarrow pC^*(E,\CL,\CE )p $$
such that $\pi(t_i)=s_{\af_i}p_{A_i}$ for $1 \leq i \leq n$. 
It is then immediate to have 
$\pi(\bt_z(t_i))=\gm_z(\pi(t_i))$ for all $i$ 
and the gauge invariant theorem 
(\cite[Theorem 5.3]{BP1}) proves the injectivity of $\pi$.
\end{proof}

\vskip 1pc

Several equivalent conditions for a graph $C^*$-algebra 
$C^*(E)$ to be purely infinite are known in \cite[Theorem 2.3]{HZ}. 
In particular, if $E$ is a row-finite graph, then 
$C^*(E)$ is purely infinite if and only if 
$E$ satisfies condition (K) and every vertex 
in a maximal tail connects to a loop in the maximal tail. 
Thus if $\af$ is a loop in $E$, there is 
another loop $\bt$ which is neither an initial path of $\af$ 
nor an extension of $\af$ whenever 
$C^*(E)$ is purely infinite. 
In the same vein, we have the following proposition.  
 
\vskip 1pc
  
\begin{prop}\label{ml} 
Let $(E,\CL,\CE)$ be a labeled space with a loop 
at a minimal set $A\in \CE_{\rm min}$. 
If $C^*(E,\CL,\CE)$ is purely infinite, then  
there exist at least two loops $(\af,A)$ and $(\bt,A)$ such that 
$\bt^m\neq \af^k$ for all $m,k\geq 1$.
\end{prop}

\begin{proof} 
Let $(\af,A)$ have the smallest length among the loops at $A$. 
To prove the theorem, we  show that 
if there is no loops at $A$ other than $(\af^k,A)$,  
then a quotient algebra of 
$C^*(E,\CL,\CE)/I_H$ (for a hereditary saturated subset $H$ of $\CE$) 
has a hereditary subalgebra isomorphic to $M_n(C(\mathbb T))$, 
which then contradicts to the assumption that 
$C^*(E,\CL,\CE)$ is purely infinite. 

Observe   for the minimal set $A\in \CE_{\rm min}$ that
$$(\bt,A) \text{ is a loop if and only if }  
\, \bt\in \CL(AE^{\geq 1}A),$$ 
which is immediate from 
$A\subset r(A,\bt)$ if and only if $r(A,\bt)\cap A\neq \emptyset$.    
Suppose to the contrary that 
\begin{eqnarray}\label{oa} 
\text{if }\bt\in\CL(AE^{\geq 1}A), 
\text{ then }  \bt^m\in \,\{ \alpha ^k  :\, k \geq 1 \} 
\ \text{ for some } m\geq 1,
\end{eqnarray}
and let
$$H:= \big\{B\in \CE: r(B,\bt)\cap A=\emptyset 
\ \text{ for all } \bt\in \CL^*(E)\big\}.$$ 
Then $H\neq \CE$ because $A \notin H$. 
Since $C^*(E,\CL,\CE)$ is purely infinite, by Lemma~\ref{MLQ}, 
$\af$ must have an exit, hence one of the following holds.
\begin{enumerate}
\item[(i)] There exists a $\gm \in \CL(AE^{|\af|})$ with  $\gm\neq \af$.
\item[(ii)] $A \subsetneq r(A,\af)$.
\end{enumerate}
If (i) holds,  then $r(A,\gm)\in H$ since 
$r(r(A,\gm),\bt)\cap A=r(A,\gm\bt)\cap A\neq \emptyset$ 
for some $\bt$ implies that $\gm\bt\in \CL(AE^{\geq 1}A)$, 
that is,  $\gm\bt$ is a loop at $A$. 
Hence 
$(\gm\bt)^m=\af^k$ for some $m,k\geq 1$  by (\ref{oa}),  
which then gives $\gm=\af$, a contradiction. 
If (ii) is the case, then $r(A,\af) \setminus A  \in H$. 
For this, suppose $B:=r(A,\af)\setminus A\notin H$. 
Then  there is $\bt\in \CL^*(E)$ such that 
$r(B,\bt)\cap A\neq \emptyset$. 
Again by minimality of $A$, we have 
\begin{eqnarray}\label{ab}
A\subset r(B,\bt)= r\big(r(A,\af)\setminus A,\bt\big)\subset r(A,\af\bt)
\end{eqnarray}
which shows that $(\af\bt,A)$ is a loop, and  
\begin{eqnarray}\label{ba}
A\subset r(A,\af)\subset r\big(r(B,\bt),\af\big)
=r(B,\bt\af)=r\big(r(A,\af)\setminus A,\bt\af\big). 
\end{eqnarray}
Again by (\ref{oa}), there exist $m,k\geq 1$ such that 
$(\af\bt)^m=\af^k$, 
but then $\af(\bt\af)^{m-1}\bt=(\af\bt)^m= \af^k=\af\af^{k-1}$ 
and thus $(\bt\af)^{m-1}\bt=\af^{k-1}$. 
Multiplying  both sides by $\af$ from the right, 
we have  
$(\bt\af)^{m}=\af^k=(\af\bt)^m$, or 
$\bt\af=\af\bt$.  
Thus we have from (\ref{ab}) and (\ref{ba}) 
a contradiction to  our standing 
assumption that $(E,\CL,\CE)$ is weakly left-resolving, and  
see that $H$ is nonempty. 

 To show that $H$ is hereditary, let $B \in H$ and  $\gm \in \CL^*(E)$. 
 If $r(B, \gm) \notin H$, 
 then $r(r(B,\gm),\bt) \cap A \neq \emptyset$ for some $\bt \in \CL^*(E)$, but then 
 $r(B,\gm\bt) \cap A \neq \emptyset$, which is  a contradiction to $B\in H$. 
 Thus $H$ is closed under relative ranges.   
 It is rather obvious that $H$ is closed under taking finite unions and subsets. 
 
  In order to see that $H$ is saturated, let 
 $B \in \CE $ satisfy $r(B,\gm) \in H$ for all  $\gm \in \CL^*(E)$. 
 If $B \notin H$, then $r(B,\bt) \cap A \neq \emptyset$ 
 for some $\bt$. 
 Thus again from minimality of $A$, 
we obtain the inclusion $A\subset r(B,\bt)$   and thus
 $$A\subset r(A,\af)\subset r\big(r(B,\bt),\af\big) $$ 
which means  $r(B,\bt)\notin H$, a contradiction.

Consider the quotient labeled space $(E,\CL,\CE/H)$. 
 If $A\cup W=\emptyset \cup W$ for some $W\in  H$, 
 then $A\in H$ (since $H$ is closed under taking subsets) 
 which is not true, hence $[A] \neq [\emptyset]$ in $\CE/H$. 
Since $r(A,\af)\setminus A$ belongs to $H$ as seen above, 
$[A] = [r(A,\af)](=r([A],\af))$,  
and thus  $(\af,[A])$ is a loop in $(E,\CL,\CE/H)$ without 
an exit of  type (ii). 
If it has an exit of type (i) with a path 
$\mu\in \CL_H([A]E^{|\af|})$, $\mu\neq \af$ and 
$r([A],\mu)\neq [\emptyset]$, then 
$[r(A,\mu)]\neq [\emptyset]$, that is, $r(A,\mu)\notin H$. 
Then  there exists a $\bt\in \CL^*(E)$ such that
$r(r(A,\mu),\bt)\cap  A\neq \emptyset$, or equivalently 
$A\subset r(A,\mu\bt)$. 
But then by (\ref{oa})  
$(\mu\bt)^m=\af^n$ for some $m,n\geq 1$.  
Thus we have $\mu=\af$, a contradiction.
 
Obviously  $[A]\in \CE/H$ is a minimal set, 
one can apply the same arguments  in the proof of  Lemma~\ref{MLQ} to 
see that the  $C^*$-algebra $C^*(E, \CL, \CE/H)$ 
contains a hereditary $C^*$-subalgebra isomorphic to 
$M_n(C(\mathbb{T}))$ for some $n\geq 1$. 
Then the quotient $C^*$-algebra 
$C^*(E,\CL,\CE)/I_H (\cong C^*(E, \CL, \CE/H))$  
is not purely infinite, a contradiction.
\end{proof}

\vskip 1pc 
 
 Since pure infiniteness of $C^*(E)$ implies that
every quotient graph of $E$ satisfies condition (L), 
one might expect  that 
a labeled space $(E,\CL,\CE)$ should be strongly disagreeable  
whenever its $C^*$-algebra $C^*(E,\CL,\CE)$ is purely infinite. 
In fact, this is true in some cases 
including all (locally finite directed) graphs.

\vskip 1pc

\begin{prop}\label{prop-nondisagreeable} 
Let $(E,\CL,\CE)$ be a labeled space such that 
every quotient labeled space satisfies the condition $(*)$.
If $C^*(E,\CL,\CE)$ is purely infinite, 
then $(E,\CL,\CE)$ is strongly disagreeable.
\end{prop}

\begin{proof} 
We prove the proposition for a labeled space $(E,\CL,\CE)$  
since the argument we give here is valid for  quotient labeled spaces. 

Suppose $(E,\CL,\CE)$ is not disagreeable. 
By Proposition~\ref{prop-disagreeable} (or Remark~\ref{q-disagreeable}), 
there exist $A\in\CE$ and $\dt\in \CL^*(E)$ 
satisfying $\CL(AE^{|\dt|n})=\{\dt^n\}$ for all $n\geq 1$. 
By choosing $\dt$ with the smallest length, 
this is equivalent to say that 
\begin{eqnarray}\label{path from A} 
\CL(AE^{\geq 1})=\{\dt^n\dt': 
  n\geq 0,\ \dt' \text{ is an initial path of } \dt\}. 
\end{eqnarray} 
Since $A$ contains a minimal subset and  
(\ref{path from A}) holds true for any subset of $A$ in $\CE$,  
we may assume that $A$ is minimal. 
Moreover, any loop at $A$ is of the form $\dt^n$. 
In fact, if $\af= \dt^n \dt'$ is a loop 
for some $n \geq 1$ and $\dt'$ with $\dt=\dt'\dt''$,   
then $\af^{\infty}=(\dt^n\dt')^{\infty}=\dt^{\infty}$ by (\ref{path from A}), 
which implies $\dt'\dt''=\dt''\dt'$.  
Since $\dt$ has the smallest length 
among the paths satisfying $(\ref{path from A})$, 
we must have  $\dt=\dt'$   (see the proof of Lemma \ref{path comparison}).

We claim  that  $\CL(AE^{\geq 1}A)= \emptyset$. 
If there is a path $\af\in \CL(AE^{\geq 1}A)$, then  
$\af$ satisfies $r(A,\af)\cap A\neq \emptyset$, and  hence 
it is  a loop at $A$ since $A$ is minimal. 
We may assume $\af$ has the smallest length among the paths 
in  $\CL(AE^{\geq 1}A)$.   
Since $C^*(E,\CL,\CE)$ is purely infinite, 
by Proposition~\ref{ml}  there exists a loop $\bt$ at $A$ 
such that $\bt^m \notin \{\alpha^k: k\geq 1\}$ for all $m \geq 1$.
Thus   $\bt \neq \dt^m$ for all $m \geq 1$,  a contradiction.

Let $\CE_{\rm min}$ be the set of all minimal sets in $\CE$ and 
let 
$$\{B_n : n\geq 1\}:= 
\{B\in \CE_{\rm min}: B \subset r(A,\af) \text{ for some } 
\af\in \CL(AE^{\geq 1}) \}.$$ 
as in the proof of Theorem~\ref{thm-IH}. 
Then the previous arguments for $A$ can be 
applied to $B_n$ to see  that 
\begin{eqnarray}\label{Bn}
r(B_n,\af)\cap B_n=\emptyset
\end{eqnarray} 
for all $n\geq 1$ and $\af\in \CL(B_nE^{\geq 1})$. 
If $n\neq m$ and $B_n\cap r(B_m, \ld) \neq \emptyset$ for some $\gm\in \CL^*(E)$ 
(hence, if  $B_n\subseteq r(B_m, \ld)$), 
we write $B_m>B_n$. 
Note that once $B_m>B_n$ holds, 
$B_n>B_m$ is not possible because of (\ref{Bn}).
Then the hereditary subset $H$   of $\CE$ consisting of 
all finite unions of $B_n$'s gives rise to an ideal 
$I_H$, and this ideal can be shown to be an AF algebra 
in the same way as before, 
which is a contradiction to Remark~\ref{rmk-pre}.(c). 
\end{proof}

\vskip 1pc

In view of the fact that if $C^*(E,\CL,\CE)$ is purely infinite, then 
every quotient $C^*$-algebra (particularly
every $C^*$-algebra of a quotient labeled space)  
is purely infinite,  
Proposition~\ref{prop-nondisagreeable} supports our stance to discuss 
the pure infiniteness of labeled graph $C^*$-algebras only 
for strongly disagreeable labeled spaces. 
Theorem~\ref{SDPI} below is a non-simple version of 
Proposition~\ref{simple PI}. 
Note that the equivalent conditions in 
the theorem are also equivalent to the following: 
every nonzero hereditary $C^*$-subalgebra of 
$C^*(E,\CL,\CE/H )$ contains an infinite projection 
for every hereditary saturated subset $H$ of $\CE$. 

\vskip 1pc
 
\begin{thm}\label{SDPI} 
Let $(E, \CL, \CE )$ be a strongly disagreeable labeled space. 
Then the following are equivalent: 
\begin{enumerate}
\item[(a)] Every nonzero hereditary $C^*$-subalgebra in every quotient of 
$C^*(E,\CL,\CE )$ contains an infinite projection. 
\item[(b)]  $C^*(E,\CL,\CE )$  is purely infinite.
\item[(c)] For each $A \in \CE $ the projection $p_A$ is properly infinite. 
\end{enumerate}
\end{thm}

\vskip 1pc
\begin{proof} We only need to show that (c) implies (a).  
If $C^*(E,\CL,\CE )$ is simple, the assertion follows from 
Proposition~\ref{simple PI}. 
If $C^*(E,\CL,\CE )=C^*(s_a,p_A)$ is not simple, and 
$I$ is a nonzero ideal of  $C^*(E,\CL,\CE )$, 
then by Lemma \ref{GI}, $I$ is gauge-invariant, that is, $I=I_H$ 
for a nonempty hereditary saturated subset $H$ 
(in fact, $H=\{A \in \CE :p_A \in I\}$).
Moreover, with $C^*(E,\CL, \CE/H)=C^*(t_a, p_{[A]})$, 
there is an isomorphism 
$$\psi:  C^*(E,\CL, \CE/H)\to C^*(E,\CL,\CE)/I$$ 
such that $\psi(p_{[A]})=p_A+I$ and $\psi(t_a)=s_a+I$ for 
$[A]\in \CE/H$ and $a\in \CA_H$ (Remark~\ref{remark-review}). 
 Now, we show that a nonzero hereditary $C^*$-subalgebra  
$B$ of $C^*(E,\CL, \CE/H)$ has an infinite projection. 
 Since  $(E,\CL, \CE/H)$ is disagreeable, 
applying Proposition~\ref{DHP} 
 (also see Remarks~\ref{rmk of DHP}.(a)) 
 we have  a nonzero projection $p\in B$ such that 
 $s_{\mu}p_{[A]}s_{\mu}^* \preceq p$ for some 
 $\mu \in \CL_H^*(E)$ and 
 $[A]\in  \CE/H$, $[A]\neq [\emptyset]$. 
 On the other hand, 
 the projection $p_A +I \neq 0$ is properly infinite 
 in $C^*(E,\CL,\CE)/I$  
 because $p_A$ is properly infinite by (c). 
But the property of being properly infinite is preserved 
under a $*$-homomorphism, 
hence the projection $p_{[A]}=\psi^{-1}(p_A+I)$ is properly infinite. 
Thus  $p$ is infinite in $B$. 
\end{proof}

\vskip 1pc

\begin{cor}\label{HCL} 
Let $(E,\CL,\CE )$ be a strongly disagreeable labeled space. 
If for every hereditary saturated subset $H$ of $\CE $,
every vertex connects to a loop 
in the quotient labeled space $(E,\CL,\CE/H)$, 
then $C^*(E,\CL,\CE)$ is purely infinite. 
\end{cor}
\begin{proof} 
This follows from Theorem~\ref{SDPI}  and Remarks~\ref{rmk of DHP} .
\end{proof}
  
\vskip 1pc

\begin{ex} 
For the labeled graph $(E,\CL)$  given by 
\vskip 1pc
\hskip 4pc
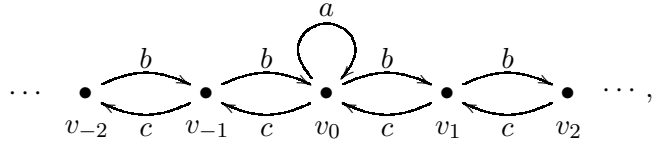
\xy
/r0.38pc/:
(-25,0)*+{\cdots};(25,0)*+{\cdots,};
(-20,0)*+{\bullet}="V-2";
(-10,0)*+{\bullet}="V-1";
(0,0)*+{\bullet}="V0";
(10,0)*+{\bullet}="V1";
(20,0)*+{\bullet}="V2";
"V-2";"V-1"**\crv{(-20,0)&(-15,3.5)&(-10,0)};
?>*\dir{>}\POS?(.5)*+!D{};
"V-1";"V0"**\crv{(-10,0)&(-5,3.5)&(0,0)};
?>*\dir{>}\POS?(.5)*+!D{};
"V0";"V1"**\crv{(0,0)&(5,3.5)&(10,0)};
?>*\dir{>}\POS?(.5)*+!D{};
"V1";"V2"**\crv{(10,0)&(15,3.5)&(20,0)};
?>*\dir{>}\POS?(.5)*+!D{};
"V-1";"V-2"**\crv{(-10,0)&(-15,-3.5)&(-20,0)};
?>*\dir{>}\POS?(.5)*+!D{};
"V0";"V-1"**\crv{(0,0)&(-5,-3.5)&(-10,0)};
?>*\dir{>}\POS?(.5)*+!D{};
"V1";"V0"**\crv{(10,0)&(5,-3.5)&(0,0)};
?>*\dir{>}\POS?(.5)*+!D{};
"V2";"V1"**\crv{(20,0)&(15,-3.5)&(10,0)};
?>*\dir{>}\POS?(.5)*+!D{};
"V0";"V0"**\crv{(0,0)&(-4,3.5)&(0,7)&(4,3.5)&(0,0)};
?>*\dir{>}\POS?(.5)*+!D{};
(-15,3)*+{b};(-5,3)*+{b};(5,3)*+{b};(15,3)*+{b};
(-15,-3)*+{c};(-5,-3)*+{c};(5,-3)*+{c};(15,-3)*+{c};
(0,7)*+{a};(0.1,-3)*+{v_0};(10.1,-3)*+{v_1};
(-9.9,-3)*+{v_{-1}};
(-19.9,-3)*+{v_{-2}};
(20.1,-3)*+{v_{2}};
\endxy
\vskip 1pc
\noindent  
the accommodating set $\CE$ consists of the finite or co-finite 
vertex subsets, and 
 $$H: =\{A \subset E^0  : A ~\text{is finite}\}$$ is the 
only non-trivial hereditary and saturated subset of $\CE$.
Each single vertex set $\{v_i\}=[v_i]_{|i|} \in \CE $, 
$i\in \mathbb Z$ admits a loop, hence the projections 
$p_{\{v_i\}}$ are all infinite because
the labeled space $(E,\CL,\CE )$ is disagreeable.   
Since $\CA_H=\{b,c\}$, $\CE/H=\{[\emptyset], [E^0]\}$, and 
$[r(b)]=[r(c)]=[E^0]$, 
the quotient labeled space $(E,\CL, \CE/H)$ can be 
visualized as follows.
\vskip 1pc
\hskip 10.5pc
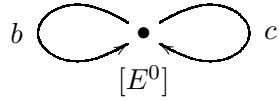
\xy
/r0.5pc/:
(0,0)*+{\bullet}="V0";
"V0";"V0"**\crv{(0,0)&(4,3)&(8,0)&(4,-3)&(0,0)};?>*\dir{>}\POS?(.5)*+!D{};
"V0";"V0"**\crv{(0,0)&(-4,3)&(-8,0)&(-4,-3)&(0,0)};?>*\dir{>}\POS?(.5)*+!D{};
(-8,0)*+{b},(8,0)*+{c},(0,-3)*+{[E^0]},
\endxy
\vskip 1pc
\noindent
Clearly $(E,\CL,\CE/H)$ is disagreeable,  
hence the labeled space $(E,\CL,\CE)$ is strongly disagreeable.
For each $A\in \CE$, the infinite projection $p_A\in C^*(E,\CL,\CE)$ 
is either zero or infinite in the quotient algebra 
$C^*(E,\CL,\CE)/I_H\cong C^*(E,\CL,\CE/H)$, and thus  
$p_A$ is properly infinite (\cite[Proposition 3.14]{KR}). 
Therefore the $C^*$-algebra $C^*(E,\CL,\CE)$ is 
purely infinite by Theorem~\ref{SDPI}.
\end{ex}

\end{document}